\documentclass[10pt,a4paper]{article}
\usepackage{amsfonts}
\usepackage{amsmath}
\usepackage{graphicx}
\usepackage{latexsym,amssymb}
\usepackage{amsmath}
\newtheorem{theorem}{Theorem}[section]

\newtheorem{remark}[theorem]{Remark}

\newtheorem{corollary}[theorem]{Corollary}
\newtheorem{definition}[theorem]{Definition}

\newtheorem{lemma}[theorem]{Lemma}

\newtheorem{proposition}[theorem]{Proposition}
\newenvironment{proof}[1][Proof]{\emph{#1.} }{\hfill{\mbox{$\square$}}}
\begin{document}
  \title{\textbf{Lightlike hypersurfaces in indefinite $\mathcal{S}$-manifolds}}
 \author{\begin{tabular}{cc}
  Letizia Brunetti and Anna Maria Pastore
\end{tabular}}
\date{}
\maketitle

\begin{abstract}
\noindent In a metric $g.f.f$-manifold we study lightlike hypersurfaces $M$ tangent to the characteristic vector fields, and owing to the presence of the $f$-structure, we determine some decompositions of $TM$ and of a chosen screen distribution obtaining two distributions invariant with respect to the structure. 
We discuss the existence of a $g.f.f$-structure on a lightlike hypersurface and, under suitable hypotheses, we obtain an indefinite $\mathcal{S}$-structure on the leaves of an integrable distribution. The existence of totally umbilical lightlike hypersurfaces of an indefinite $\mathcal{S}$-space form is also discussed. Finally, we explicitely describe a lightlike hypersurface of an indefinite $\mathcal{S}$-manifold.
\end{abstract}
\textbf{2000 Mathematics Subject Classification} 53C40, 53C50, 53D10
\\
\textbf{Keywords and phrases} Lightlike hypersurfaces, indefinite globally framed $f$-structures

\section{Introduction}
In the wide field of semi-Riemannian geometry, the study of
lightlike, or degenerate, submanifolds comes now to fill
a gap in the general theory of submanifolds. In fact, while the geometry
of submanifolds of Riemannian manifolds, since \cite{Ch1},
 has received powerful impulse, polarizing a lot of attention,
with studies conducted in great generality and developing a great variety
of techniques, on the contrary the study of degenerate geometry is a
relatively new field of research. It rises within the semi-Riemannian
context, due to the existence of the so-called \emph{causal character} of
geometrical objects: their spacelike, timelike or lightlike nature, in
fact, implies the existence of three types of hypersurfaces and
submanifolds. Among them, the spacelike and timelike cases have received a
systematic exposition in the fundamental book \cite{[O'N]}. We have to be
looking forward 1990 to find the first studies about the lightlike case,
when A.\ Bejancu and K.L.\ Duggal introduced the lightlike geometry
(\cite{Bj1,Bj2,DB1,BF,BFL,DB}). At the moment lightlike hypersurfaces are studied in paraquaternionic K\"ahler manifolds'contexts by S. Ianu\c s, R. Mazzocco\ and\ G. E. V\^ilcu (\cite{IMV}). 

The primary difference between the lightlike submanifolds and
non-degenerate submanifolds arises due to the fact that in the first case
the normal vector bundle has non trivial intersection with the tangent vector bundle, and moreover in a lightlike hypersurface the normal vector bundle is contained in the tangent vector bundle. Thus,
one fails to use the classical theory of non-degenerate submanifolds to
define the induced geometrical objects (as linear connection, second
fundamental form, Gauss and Weingarten equations) on a lightlike
hypersurface.

The growing importance of lightlike geometry is 
motivated by its extensive use in mathematical physics, in particular in relativity.
In fact, semi-Riemannian manifolds $(\bar{M},\bar{g})$ with $dim \bar{M}=n>4$ are 
natural generalizations of spacetime of general relativity and lightlike hypersurfaces 
are models of different types of horizons separating domains of $(\bar{M},\bar{g})$ with different 
physical properties.
 
There exist also reasons that motivate the study of the lightlike hypersurfaces of indefinite $g.f.f$-manifolds, in particular of indefinite $\mathcal{S}$-manifolds. In \cite{DB}, K.L.\ Duggal and A.\ Bejancu proved that a lightlike framed hypersurface of a Lorentz $\mathcal{C}$-manifold, with an induced metric connection, is a Killing horizon. In a recent paper (\cite{DS}), K.L.\ Duggal and B.\ Sahin begin to work on lightlike submanifolds of indefinite Sasakian manifolds because the contact geometry has a significant use in differential equations, optics and phase spaces of dynamical systems. Furthermore, in \cite{Dug0}, K.L.\ Duggal shows that a globally hyperbolic spacetime and the de Sitter spacetime can carry a framed structure.
  
We begin with some basic information about indefinite $\mathcal{S}$-manifolds and about lightlike hypersurfaces of a semi-Riemannian manifold. Afterwards, for a metric $g.f.f$-manifold we consider a lightlike hypersurface $M$ tangent to the characteristic vector fields, we introduce a particular screen distribution $S(TM)$, using the properties of the indefinite $\mathcal{S}$-manifold. Then, finding other decompositions of $S(TM)$ and $TM$ yields two distributions $D_{0}$ and $D$ on $M$ that are studied in section~4. 
We discuss the existence of a $g.f.f$-structure on a lightlike hypersurface and we obtain an indefinite $\mathcal{S}$-structure on the leaves of $D_{0}$, if $D_{0}$ is an integrable distribution. Section 5 deals with the existence of totally umbilical lightlike hypersurfaces of an indefinite $\mathcal{S}$-space form.
In the last section, we consider the three examples of indefinite $\mathcal{S}$-manifolds given in \cite{LP}. For the first one we  explicitely describe a lightlike hypersurface, to which we apply the previous results, while, for the other two examples, we prove that they do not admit lightlike hypersurfaces tangent to the characteristic vector fields.

All manifolds, tensor fields and maps are assumed to be smooth. We shall use the Einstein convention omitting the sum symbol for repeated indexes.

\emph{Acknowledgments}. The authors are grateful to Prof. S. Ianus for discussions about the topic of this paper during his stay at the University of Bari and the stay of the first author at the University of Bucharest.

\section{Preliminaries}
A manifold $\bar{M}$ is called a \emph{globally framed $f$-manifold} if it is endowed with a non null $(1,1)$-tensor field $\bar{\varphi}$ of constant rank, such that $\ker\bar{\varphi}$ is parallelizable i.e. there exist global vector fields  $\bar{\xi}_{\alpha}$, $\alpha\in\{1,\ldots,r\}$, with their dual $1$-forms $\bar{\eta}^{\alpha}$, satisfying $\bar{\varphi}^{2}=-I+\bar{\eta}^{\alpha}\otimes\bar{\xi}_{\alpha}$ and $\bar{\eta}^{\alpha}(\bar{\xi}_{\beta})=\delta_{\beta}^{\alpha}$.

The $g.f.f$-manifold $(\bar{M}^{2n+r},\bar{\varphi},\bar{\xi}_{\alpha}
,\bar{\eta}^{\alpha}, \bar g) $, $\alpha \in \{1,\ldots,r\}$ is said an indefinite metric $g.f.f$-manifold if $\bar g$ is a semi-Riemannian metric, with index $\nu$, $0<\nu < 2n+r$, satisfying the following compatibility condition
 \begin{equation*}
\bar{g}(\bar{\varphi}X,\bar{\varphi}Y)=\bar{g}(
X,Y)-\varepsilon_{\alpha}\bar{\eta}^{\alpha}(
X)  \bar{\eta}^{\alpha}(Y)
\end{equation*}
for any $X,Y\in \Gamma T\bar{M}$,  being $\varepsilon_{\alpha}=\pm1$ according to whether $\bar{\xi}_{\alpha}$ is
spacelike or timelike. Then, for any
$\alpha\in\{1,\ldots,r\}$, one has $\bar{\eta}^{\alpha}(X)=\varepsilon_{\alpha}\bar{g}(X,\bar{\xi}_{\alpha})$.

An indefinite metric $g.f.f$-manifold is called
 \emph{indefinite $\mathcal{S}$-manifold} if it is normal and $d\bar{\eta}^{\alpha}=\Phi$, for any $\alpha\in\{1,\ldots,r\}$,
where $\Phi(X,Y)  =\bar{g}(X,\bar{\varphi}Y)$ for any $X,Y\in \Gamma(T\bar{M})$.
 The normality condition is expressed by the vanishing of the tensor field $N=N_{\bar\varphi} +2d\bar\eta^{\alpha}\otimes \bar\xi_{\alpha}$, $N_{\bar\varphi}$ being the Nijenhuis torsion of $\bar\varphi$.
Furthermore, as proved in \cite{LP}, the Levi-Civita connection of an indefinite $\mathcal{S}$-manifold
satisfies:
 \begin{equation}\label{LC}
(  {\bar\nabla}_{X}{\bar\varphi})  Y={\bar g}(  {\bar\varphi
}X,{\bar\varphi}Y) \bar{\xi}+\bar{\eta}(  Y)
{\bar\varphi}^{2}(  X).
\end{equation}
 where $\bar\xi = \sum_{\alpha =1}^r \bar\xi_{\alpha}$  and $\bar\eta =\sum_{\alpha=1}^r \varepsilon_{\alpha}\eta^{\alpha}$.
We recall that $\bar\nabla_X \bar\xi_{\alpha}=-\varepsilon_{\alpha}\bar\varphi X$ and $\ker \bar\varphi $ is a integrable flat distribution since $\bar\nabla_{\bar\xi_{\alpha}} \bar\xi_{\beta} =0$. For more details we refer to \cite{LP}.

\bigskip
Fol\-lowing \cite{DB}, we re\-call some ba\-sic re\-sults about lightlike hyper\-sur\-faces  of a semi-Riemannian ma\-ni\-fold $(\bar{M},\bar{g})$. Given a lightlike hypersurface $M$ of $\bar{M}$, one can consider for any $p\in M$
the vector spaces:
 \begin{align*}
T_{p}M^{\bot}&=\{ E_p\in T_{p}\bar{M}\,\,\mid\bar{g}_{p}(E_{p},W)=0\,\,\,\,\,\textrm{for all} \,\, W\in T_{p}M\},\\
Rad(T_{p}M) &= \{  V\in T_{p}M\,\mid\,g_{p}(  V,W)
=0\,\,\,\,\,\textrm{for all} \,\, W\in T_{p}M\}  =T_{p}M^{\bot}\cap T_{p}M.
\end{align*}
 Then, $Rad(T_{p}M)=T_{p}M^{\bot}\subset T_{p}M$ and one has the $1$-dimensional degenerate distribution on $M$, called the {\it radical distribution}. A {\it screen distribution} on $M$ is defined as a distribution complementary to the radical one, so that we have
\begin{align}
TM& = Rad(TM)\bot S(TM)\label{decomp01}\\
 T\bar{M}|_{M}& = S(TM)\bot S(TM)^{\bot}\nonumber
\end{align}
where $S(TM)^{\bot}$ is the complementary vector bundle
to $S(  TM)  $ in $ T\bar{M}|  _{M}$ with respect to
$\bar{g}$.
Obviously, there exist several screen distributions and they are non-degenerate.

We report the following theorem proved in \cite{DB}, adapting it to our context. 
\begin{theorem}[\cite{DB}]
\label{teorsezioneN}Let $(
M,g,S(  TM)  )  $ be a lightlike hypersurface of an indefinite $g.f.f$-manifold $(  \bar{M},\bar{\varphi},\bar{\xi}_{\alpha},\bar{\eta}^{\alpha},\bar{g})  $. Then there exists
a unique rank one vector subbundle $ltr(  M)  $ of $T\bar{M}$, with
base space $M$, such that for any non-zero section $E$ of $Rad(TM)$ on a
coordinate neighbourhood $\mathcal{U}\subset M$, there exists a unique section
$N$ of $ltr(  M)  $ on $\mathcal{U}$ satisfying:
\begin{equation}
\bar{g}(  N,E)  =1,\quad \bar{g}(  N,N)  =0,\qquad \bar{g}(  N,W)  =0\,\,\textrm{  for all }\, W\in\Gamma(  S(  TM)  |  _{\mathcal{U}
}). \label{equ04}
\end{equation}
$ltr(  M)  $ is called the {\it lightlike
transversal vector bundle} of $M$ with respect to $S(  TM)  $.
\end{theorem}
One can consider the following decompositions
\begin{align}
S(  TM)  ^{\bot}&=Rad(TM)\oplus ltr(M)
,\label{decomp03}\\
T\bar{M}|  _{M}&=TM\oplus ltr(  M).\label{decomp04}
\end{align}
 Let $\bar{\nabla}$ be the Levi-Civita connection on $\bar{M}$. Using (\ref{decomp04}) we deduce 
\begin{equation*}
\bar{\nabla}_{X}Y=\nabla_{X}Y+h(  X,Y), \qquad \bar{\nabla}_{X}V=-A_{V}X+\nabla_{X}^{lt}V,
\end{equation*}
for any $X,Y\in\Gamma(  TM)  $ and $V\in\Gamma(  ltrM)
$.
\noindent Following \cite{DB}, $\nabla$ and $\nabla^{lt}$ are called the \emph{induced connections }on $M$ and
$ltr(  M)  $ respectively, and as in the classical theory of
Riemannian hypersurfaces, $h$ and $A_{V}$ are called the \emph{second
fundamental form }and the \emph{shape operator}, respectively. Further, the above equations are cited as the \emph{Gauss }and \emph{Weingarten
equation,} respectively. Locally, let $E$, $N$ and $\mathcal{U}$ be as in
Theorem \ref{teorsezioneN}, then for any $X,Y\in\Gamma( 
TM|  _{\mathcal{U}})  $, putting:
\begin{equation}
B(  X,Y)  =\bar{g}(  h(  X,Y)  ,E), \quad {\rm and}\quad \tau(  X)  =\bar{g}(  \nabla_{X}^{lt}N,E)
,\nonumber%
\end{equation}
for any $X,Y\in\Gamma(TM|_{\mathcal{U}})  $, we can write:
\begin{equation}
\bar{\nabla}_{X}Y=\nabla_{X}Y+B(  X,Y)  N,\quad {\rm and} \quad \bar{\nabla}_{X}N=-A_{N}X+\tau(  X)  N,\label{equ12}%
\end{equation}
$B$ is called the \emph{local second fundamental form }of $M$,
because it determines $h$ on $\mathcal{U}$. 
As proved in \cite{DB1}, the local second fundamental form of $M$ on $\mathcal U$ is independent of the choice of the screen distribution. Moreover $B$ is degenerate and $B(X,E)=0$ for any $X \in \Gamma(TM|_{\mathcal{U}})$.

The decomposition (\ref{decomp01}) allows to define a canonical projection 
$P:\Gamma(TM)\rightarrow\Gamma(S(
TM)  )$.

\noindent Then, for any $X,Y\in\Gamma(TM)$ and
$U\in\Gamma(  Rad(TM))  $ we can write
\begin{equation}
\nabla_{X}PY=\overset{*}{\nabla}_{X}PY+\overset{*}{h}(  X,PY)
, \quad {\rm and} \quad \nabla_{X}U=-\overset{*}{A}_{U}X+\overset{*}{\nabla_{X}^{t}}U,\label{equ14}%
\end{equation}
where $\overset{\ast}{\nabla}$ and $\overset{\ast}{\nabla^{t}}$ are linear connections on the bundles $S(TM)$ and $Rad(TM)$, respectively. Further, $\overset{\ast}{h}:\Gamma(TM)\times\Gamma(S(TM))\rightarrow\Gamma(Rad(TM))$ is $\mathfrak{F}(M)$-bilinear and $\overset{\ast}{A}_{U}:\Gamma(TM) \rightarrow\Gamma(S(TM))$ is an 
$\mathfrak{F}(M)$-linear operator and they are called the \emph{second
fundamental form }and the \emph{shape operator }of the screen distribution, respectively. The equations in (\ref{equ14}) are cited as the \emph{Gauss equation }and the \emph{Weingarten equation}. 
Locally, let $\mathcal{U}$ be a coordinate neighbourhood of $M$, and $E$, $N$
sections on $\mathcal{U}$, as in Theorem \ref{teorsezioneN}. Then, putting
$C(  X,PY)  =\bar{g}(  \overset{*}{h}(  X,PY)
,N)$ for any $X,Y\in\Gamma( TM|_{\mathcal{U}})  $, one has $\,\overset{\ast}{h}(  X,PY)  =C(  X,PY)  E$, $\bar g (\overset{*}{\nabla_{X}^{t}}E,N) =-\tau(X)$ and, locally on $\mathcal{U}$, (\ref{equ14}) becomes
\begin{equation}
\nabla_{X}PY=\overset{*}{\nabla}_{X}PY+C(  X,PY)  E, \quad {\rm and} \quad
\nabla_{X}E=-\overset{*}{A}_{E}X-\tau(  X)  E.\label{equ21}%
\end{equation}
Finally, geometrical objects of the lightlike hypersurface and of the screen distribution are related as follows, for any $X,Y\in\Gamma(TM)$:
\begin{equation}\label{num05}
g(A_{N}X,PY)=C(X,PY), \;\;  g(  \overset{*}{A}_{E}X,PY)  =B(  X,PY),\;\; \bar{g}(  A_{N}X,N)  =0,\;\; \bar{g}(  \overset{*}{A}_{E}X,N)  =0.
\end{equation}
Furthermore, one has:
$\overset{*}{A}_{E}E=0$, $\bar{\nabla}_{E}E=\nabla_{E}E=-\tau(  E)  E.$

\begin{remark}\label{remNab}
\emph{The induced connection $\nabla$ on $M$ does not depend on the choice of 
$S(  TM)  $ if and only if the second fundamental form $h$ of $M$
vanishes identically. Furthermore, $\nabla$ is not a metric connection, in
fact it satisfies
$(  \nabla_{X}g)  (  Y,Z)  =B(  X,Y)
\bar{g}(Z,N)  +B(  X,Z)  \bar{g}(Y,N)$,
for any $X,Y,Z\in\Gamma( TM|  _{\mathcal{U}})  $.
However, if we choose $Y,Z\in\Gamma(  S(
TM)  )  $, we get $(  \nabla_{X}g)  (  Y,Z)
=0$, and using (\ref{equ14}) we easily obtain that the linear connection
$\overset{*}{\nabla}$ on $S(  TM)  $ is a metric connection.
Finally, $S(TM)$ is an integrable distribution if and only if $\overset{*}{h}$ is symmetric on $S(TM)$. }
\end{remark}
\section{Characteristic lightlike hypersurfaces of indefinite $g.f.f$-manifolds}
Let $(  \bar{M},\bar{\varphi},\bar{\xi}_{\alpha},\bar{\eta}^{\alpha}%
,\bar{g})  $ be an indefinite $g.f.f$-manifold and $M$ a
lightlike hypersurface. From the existence of the $f$-structure, for any $Z \in \Gamma(Rad(TM))$, one has $\bar{g}(  \bar{\varphi}Z,Z)  =0$, therefore
$\bar{\varphi}Z \in~\Gamma(TM)$, and we get a $1$-dimensional distribution $\bar{\varphi}( Rad(TM))  $ on $M$.

Moreover, it is easy to state the following result.
 \begin{proposition}
Let $M$ be a lightlike hypersurface of an indefinite $g.f.f$-manifold $(\bar{M},\bar{\varphi},\bar{\xi}_{\alpha},\bar{\eta}^{\alpha}
,\bar{g})$ such that the characteristic vector fields $\bar{\xi}_{\alpha}$ are tangent to $M$. Let $E$ be a non zero section of $Rad(TM)$. Then there exists a screen distribution
such that $\ker\,\bar\varphi \subset S(TM)$ and $\bar{\varphi}(E)$ belongs to $\Gamma(S(TM))$. 
\end{proposition}
\begin{definition} \emph{A lightlike hypersurface $M$ of $\bar{M}$ is called \emph{characteristic} if all the characteristic vector fields $\bar{\xi}_{\alpha}$ are tangent to $M$. A screen distribution $S(TM)$ will be called \emph{characteristic} if $\ker\,\bar\varphi \subset S(TM)$ and $\bar{\varphi}(E) \in\Gamma(S(TM))$.} 
\end{definition}

From now on, we shall write simply $(M,g,S(TM))$ to denote a characteristic lightlike hypersurface $(M,g)$, together with the choices of a fixed non zero section $E$ of $Rad(TM)$, a fixed characteristic screen distribution $S(TM)$ and the $ltr(M)$ and $N$ as in Theorem \ref{teorsezioneN}. 
In these hypotheses $\bar{\varphi}N\in \Gamma(S(TM))$. Namely, the vector field $\bar{\varphi}N\in \Gamma( T\bar{M}|_{M})$,  
and using (\ref{decomp03}), we have that $\bar{\varphi}N$ is orthogonal to
$S(TM) ^{\bot}$ since $\bar{g}(\bar{\varphi}N,E)  =-\bar{g}(N,\bar{\varphi
}E)=0$, and obviously $\bar{g}( \bar{\varphi}N,N)=0$.
Moreover, from the compatibility condition, we obtain $\bar{g}(  \bar{\varphi}N,\bar{\varphi}E)=1$. 
\begin{remark}\label{dimensione}\emph{
Given an indefinite $g.f.f$-manifold $(\bar{M},\bar{\varphi},\bar{\xi}_{\alpha},\bar{\eta}^{\alpha}
,\bar{g})$, the existence of a characteristic lightlike hypersurface implies that $dim\;\bar M =2n +r $ with $n \geq 2$,
being $E, N, \bar\varphi(E)$ linearly independent.} 
\end{remark}
Examples of characteristic hypersurfaces and characteristic screen distributions of indefinite $S$-manifolds are given in section \ref{example}.
It is easy to prove the following result (\cite{B})
\begin{proposition}
Let $(  \bar{M},\bar{\varphi},\bar{\xi}_{\alpha},\bar{\eta
}^{\alpha},\bar{g})$ be an indefinite $g.f.f$-manifold
and $(  M,g,S(  TM)  )$ a characteristic lightlike hypersurface of
$\bar{M}$. Then the rank 2 vector subbundle $\bar{\varphi}(  Rad(TM))
\oplus\bar{\varphi}(ltrM)$ of $S(TM)$ is non-degenerate.
\end{proposition}
Following \cite{DB}, being $S(TM)$ and $\bar{\varphi}(  Rad(TM))
\oplus\bar{\varphi}(ltrM)$ non-degenerate, we can define the unique non-degenerate distribution $D_{0}$
such that
\begin{equation}
S(  TM)  =(  \bar{\varphi}(  Rad(TM))  \oplus
\bar{\varphi}(  ltrM)  )  \bot D_{0}.\label{decomp05}%
\end{equation}
Then, each $\bar{\xi}_{\alpha}\in
D_{0}$, $D_{0}$ is $\bar{\varphi}$-invariant and, using (\ref{decomp01}), (\ref{decomp04}) and (\ref{decomp05}), we can write:
\begin{align}
TM&=D_{0}\bot {\mathcal F}\label{decomp06}\\
 T\bar{M}|  _{M}&=D_{0}\bot{\mathcal E} \label{decomp07}\\
 TM&=D\oplus \bar{\varphi}(  ltrM)\label{decomp08}
\end{align}
 where ${\mathcal F}:=(  \bar{\varphi}(  Rad(TM))  \oplus\bar
{\varphi}(  ltrM) )  \bot Rad(TM)$, ${\mathcal E}:=(\bar{\varphi}(  Rad(TM))  \oplus\bar{\varphi}(  ltrM)  )  \bot(
Rad(TM)\oplus ltr(  M)  )$ and  
$D:=D_{0}\bot\bar{\varphi}(  Rad(TM))  \bot Rad(TM)$,
 is a $\bar{\varphi}$-invariant distribution.

\medskip 
Now, we look for a $g.f.f$-structure on $(  M,g,S(  TM)  )$.
We consider the local lightlike vector fields $U:=-\bar{\varphi}N \in \bar{\varphi}(  ltrM)$ and $V:=-\bar{\varphi}E \in D$. 
From (\ref{decomp08}) any
$X\in\Gamma(  TM)  $ can be written as
\begin{equation*}
X=SX+QX \quad \text{ and }\quad QX=u(  X)  U,
\end{equation*}
where $S:TM\rightarrow D$ and $Q:TM\rightarrow \bar{\varphi}(  ltrM)$ are the canonical
projection maps, and $u$ is a local $1$-form on $M$ defined by $u(  X)  :=g(  X,V)$.
We note that 
\begin{equation}\label{rem01}
u(  U)  =1, \quad \forall Y\in\Gamma(  D)  \,\,\,\,u(  Y)
=0,\quad \bar{\varphi}^{2}N=-N.
\end{equation}
Then, applying $\bar{\varphi}$ to $X$, we obtain $\bar{\varphi}X=\bar{\varphi}(  SX)  +u(  X)
\bar{\varphi}U=\bar{\varphi}(  SX)  +u(  X)  N$.
For any $X\in~\Gamma(TM)$ we put $\varphi X:=\bar{\varphi}(  SX)$,
obtaining a tensor field $\varphi$ of type (1,1) on $M$. From the above
equality we get then
\begin{equation}
\bar{\varphi}X=\varphi X+u(  X)  N,\label{equ27}
\end{equation}
and applying again $\bar{\varphi}$, we have
$-X+\bar{\eta}^{\alpha}(  X)  \bar{\xi}_{\alpha}=\bar{\varphi
}\varphi X-u(  X)  U$.
We note that if $X\in\Gamma(TM)$, then $SX\in D$, 
$\varphi X=\bar{\varphi}(  SX)  \in D$,
so that $S(  \varphi X)  =\varphi X$.
Furthermore, since $\bar{\varphi}(  \varphi X)  =\bar{\varphi}(  S\varphi
X)  =\bar{\varphi}S(  \varphi X)  =\varphi
^{2}X$, we can write
$\varphi^{2}X=-X+\bar{\eta}^{\alpha}(  X)  \bar{\xi}_{\alpha
}+u(  X)  U$.
Finally, we have  $\varphi U=0$, since $U\in \bar{\varphi}(  ltrM)$, $\bar{\eta}^{\alpha}\circ\varphi=0$,
and $u(\varphi X)  =0$, for any $X\in\Gamma(TM)$.

Thus, we can state the following theorem.
\begin{theorem}
Let $(\bar{M},\bar{\varphi},\bar{\xi}_{\alpha},\bar{\eta}^{\alpha},\bar{g})$ be an
indefinite $\mathcal{S}$-manifold and consider $(  M,g,S(  TM))$ a characteristic lightlike hypersurface of $\bar{M}$ such that $E$ and $N$ are globally defined on $M$. Then $(M,\varphi,\bar{\xi}_{\alpha},U,\eta^{\alpha},u)$ is a $g.f.f$-manifold.
\end{theorem}

 For any $X,Y\in\Gamma(  TM)  $, we compute
the field $(\nabla_{X}\varphi)  Y$. Using (\ref{equ27}) and (\ref{equ12}), we get
\begin{align}
(  \bar{\nabla}_{X}\bar{\varphi})  Y  & =(  \nabla_{X}\varphi)  Y-u(  Y)  A_{N}X+B(
X,Y)  U +(  B(  X,\varphi Y)  +(  \nabla_{X}u)
Y+u(  Y)  \tau(  X)  )  N,\nonumber
\end{align}
then, from (\ref{LC}), comparing the components along $TM$ and $ltr(  M)  $, we have:
\begin{align}
(  \nabla_{X}\varphi)  Y  & =u(  Y)  A_{N}X-B(
X,Y)  U+\bar{g}(\bar{\varphi}X,\bar{\varphi}Y)  \bar{\xi}+\bar{\eta}(
Y) \bar{\varphi}^{2} X,\label{cond01}\\
(  \nabla_{X}u)  Y&=-B(  X,\varphi Y)  -u(
Y)  \tau(  X).\nonumber
\end{align}
\begin{definition}\emph{
Let $(  \bar{M},\bar{\varphi},\bar{\xi}_{\alpha}%
,\bar{\eta}^{\alpha},\bar{g})  $ be an indefinite $g.f.f$-manifold
and $(  M,g,S(  TM)  )  $ a lightlike hypersurface of
$\bar{M}$. Then $M$ is called \emph{totally geodesic lightlike hypersurface}
if any geodesic of $M$ with respect to the induced connection $\nabla$ is a
geodesic of $\bar{M}$ with respect to $\bar{\nabla}$.}
\end{definition}
In \cite{DB} it is proved that the previous
definition does not depend on the choice of a screen distribution and it is equivalent to the vanishing of the local second fundamental form $B$.
\begin{proposition}
Let $(  \bar{M},\bar{\varphi},\bar{\xi}_{\alpha},\bar{\eta
}^{\alpha},\bar{g})  $ be an indefinite $\mathcal{S}$-manifold and
$(  M,g,S(  TM)  )  $ a characteristic lightlike hypersurface of
$\bar{M}$. Then $M$ is totally geodesic if and only if for any $X\in
\Gamma(  TM)  $ and for any $Y\in\Gamma(  D)  $
\begin{align}
(  \nabla_{X}\varphi)  Y  &= \bar{g}(\bar{\varphi}X,\bar{\varphi}Y)  \bar{\xi}+\bar{\eta}(
Y) \bar{\varphi}^{2} X, \label{cond geod 01}\\
 A_{N}X &=-\varphi(  \nabla_{X}U)  -g(  X,U)
\bar{\xi}.\label{cond geod 02}
\end{align}
\end{proposition}

\begin{proof}
We suppose that $M$ is totally geodesic, that is $B(  X,Y)  =0$ for
any $X,Y\in \Gamma(TM)$. Therefore, from (\ref{cond01}), if $Y\in\Gamma(  D)  $
we have $u(Y)=0$ and $(\nabla_{X}\varphi)  Y=\bar{g}(\bar{\varphi}X,\bar{\varphi}Y)  \bar{\xi}+\bar{\eta}(
Y) \bar{\varphi}^{2} X$.
Again, writing (\ref{cond01}) for $U$, using (\ref{rem01}), we have
 $(  \nabla_{X}\varphi)  U =A_{N}X+\bar{g}(\bar{\varphi}X,\bar{\varphi}U)  \bar{\xi}+\bar{\eta}(
U) \bar{\varphi}^{2} X$,
from which we obtain $A_{N}X=-\varphi(  \nabla_{X}U)  -\bar
{g}(  X,U)  \bar{\xi}$, since $\varphi(U)=0$ and $\bar\eta(U)=0$.

Conversely, we suppose that the conditions (\ref{cond geod 01}) and
(\ref{cond geod 02}) hold and we prove that the local second fundamental form
$B$ vanishes. If $Y\in\Gamma(  TM)  $, using the decomposition
(\ref{decomp08}), there exists $\alpha\in\mathfrak{F}(
\mathcal{U})  $ such that $Y=Y_{d}+\alpha U$,
and, for any $X\in\Gamma(  TM)  $ we obtain $B(  X,Y)  =B(  X,Y_{d})  +\alpha B(  X,U)$.
Using (\ref{cond01}) and (\ref{cond geod 01}) with $Y=Y_{d}$, we
find $B(  X,Y_{d})  U=u(  Y_{d})  A_{N}X=0$, which implies
$B(X,Y_{d})  =0$.
From (\ref{cond01}), putting $Y=U$ and using (\ref{cond geod 02}), we get $B(X,U)  U=0,$ which implies
$B(  X,U)  =0$, concluding the proof.
\end{proof}
\section{Distributions on a lightlike hypersurface of an indefinite $\mathcal{S}$-manifold}
\subsection{The distribution $D_0$}
The following Lemma can be easily proved.
\begin{lemma}\label{cond08}
\label{Lemma S} Let $(  \bar{M},\bar{\varphi},\bar{\xi
}_{\alpha},\bar{\eta}^{\alpha},\bar{g})  $ be an indefinite
$\mathcal{S}$-manifold and $(  M,g,S(  TM)  )  $  a characteristic
lightlike hypersurface of $\bar{M}$. Let $\mathcal{U}\subset M$ be a
coordinate neighbourhood as fixed in Theorem \ref{teorsezioneN}, then for any $X,Y\in\Gamma(  D_{0})  $
\begin{equation}
\bar{g}(  (  \bar{\nabla}_{X}\bar{\varphi})  Y,E)
=0, \qquad \bar{g}(  (  \bar{\nabla}_{X}\bar{\varphi})  Y,N)=0.
\end{equation}
\end{lemma}
Now, referring to the decomposition (\ref{decomp06}), for any $X$ in $\Gamma(  TM)  $,  
$Y$ in $\Gamma(  D_{0})  $, we have:
\begin{equation}
\nabla_{X}Y=\overset{\circ}{\nabla}_{X}Y+\overset{\circ}{h}(  X,Y)
,\label{GaussD0}%
\end{equation}
where $\overset{\circ}{\nabla}$ is a linear
connection on the bundle $D_{0}$, and $\overset{\circ}{h}:\,\Gamma(  TM)  \times\Gamma(
D_{0})  \rightarrow\Gamma(  \mathcal{F})$ is $\mathfrak{F}(  M)$-bilinear. 

\noindent Let $\mathcal{U}\subset M$ be a coordinate neighbourhood as fixed in Theorem
\ref{teorsezioneN} and let $X,Y\in\Gamma( D_{0}|
_{\mathcal{U}})  $. Then, using (\ref{decomp06}),
(\ref{GaussD0}) can be written (locally) in the following way:
\begin{equation}
\nabla_{X}Y=\overset{\circ}{\nabla}_{X}Y+g( \nabla_{X}Y,\bar{\varphi}N)\bar
{\varphi}E+g( \nabla_{X}Y,\bar{\varphi}E)  \bar{\varphi}N+g( \nabla_{X}Y,N)E.\label{equ33}%
\end{equation}
Using Lemma \ref{cond08}, (\ref{equ12}), (\ref{num05}),
being $D_{0}$ $\bar{\varphi}$-invariant and $\bar{\nabla}\bar g =0$, we get
\begin{align*}
g(  \nabla_{X}Y,\bar{\varphi}N)   &  =-\bar{g}(  \bar{\varphi}(  \nabla
_{X}Y)  ,N) =-\bar{g}(  \bar{\varphi}(  \bar{\nabla}_{X}Y-B(
X,Y)  N)  ,N) \\
 &  =-\bar{g}(  \bar{\varphi}(  \bar{\nabla}_{X}Y)  ,N)
=\bar{g}(  (  \bar{\nabla}_{X}\bar{\varphi})  Y,N)
-\bar{g}(  \bar{\nabla}_{X}(  \bar{\varphi}Y)  ,N)\\
&  =-\bar{g}(  \bar{\nabla}_{X}(  \bar{\varphi}Y)  ,N)=\bar{g}(  \bar{\varphi}Y,\bar{\nabla}_{X}N)=-g(  A_{N}%
X,\bar{\varphi}Y)  =-C(  X,\bar{\varphi}Y).
\end{align*}
Again, using Lemma \ref{cond08}, (\ref{equ12}), (\ref{equ21}), (\ref{num05}),
$D_{0}$ being $\bar{\varphi}$-invariant and $\bar{\nabla} \bar g =0$
 \begin{align*}
g(  \nabla_{X}Y,\bar{\varphi}E)& =-\bar{g}(  \bar{\varphi}(  \nabla
_{X}Y)  ,E)=-\bar{g}(  \bar{\varphi}(  \bar{\nabla}_{X}Y)  ,E)
+B(  X,Y)  \bar{g}(  \bar{\varphi}N,E)=-\bar{g}(  \bar{\varphi}(  \bar{\nabla}_{X}Y)  ,E) \\
& =-\bar{g}(  \bar{\nabla}_{X}(  \bar{\varphi}Y)  ,E) =\bar{g}(  \bar{\varphi}Y,\nabla_{X}E+B(  X,E)  N) =\bar{g}(  \bar{\varphi}Y,\nabla_{X}E) =-B(  X,\bar{\varphi}Y).
\end{align*}
For any $X,Y\in\Gamma(D_{0}|  _{\mathcal{U}})  $, since $X,Y\in\Gamma(S(  TM)  |  _{\mathcal{U}})  $, we know that
$\nabla_{X}Y=\overset{*}{\nabla}_{X}Y+C(  X,Y)  E$, so we get
$g(  \nabla_{X}Y,N)  =C(X,Y)$. Therefore (\ref{equ33}) becomes
\begin{equation}
\nabla_{X}Y=\overset{\circ}{\nabla}_{X}Y-C(  X,\bar{\varphi}Y)
\bar{\varphi}E-B(  X,\bar{\varphi}Y)  \bar{\varphi}N+C(
X,Y)  E,\label{equ38}%
\end{equation}
and the local expression of $\overset{\circ}{h}$ is
\begin{equation}
\overset{\circ}{h}(  X,Y)  =-C(  X,\bar{\varphi}Y)
\bar{\varphi}E-B(  X,\bar{\varphi}Y)  \bar{\varphi}N+C(
X,Y)  E. \label{2°forma}
\end{equation}
\begin{theorem}\label{integrD0}
Let $(\bar{M},\bar{\varphi},\bar{\xi}_{\alpha},\bar{\eta}^{\alpha},\bar{g})  $ be
an indefinite $\mathcal{S}$-manifold and $(  M,g,S(  TM)
)  $ a characteristic lightlike hypersurface of $\bar{M}$, with the induced $g.f.f$-structure. The distribution
$D_{0}$ on $M$ is integrable if and only if for any
$X,Y\in\Gamma(  D_{0}) $:
\[C(X,Y)=C(Y,X),\quad C(X,\bar{\varphi}Y)=C(\bar{\varphi}X,Y),\quad B(X,\bar{\varphi}Y)=B(\bar{\varphi}X,Y).
\]
\end{theorem}
\begin{proof}
First of all, $\nabla$ being a torsion-free connection, using (\ref{equ38}), for any $X,Y\in \Gamma(D_{0})$ we get
\begin{align*}
	[X,Y]&=\overset{\circ}{\nabla}_{X}Y-\overset{\circ}{\nabla}_{Y}X + (C(Y,\bar{\varphi}X)-C(X,\bar{\varphi}Y))\bar{\varphi}E\\
	& \quad + (B(Y,\bar{\varphi}X)-B(X,\bar{\varphi}Y))\bar{\varphi}N -(C(Y,X)-C(X,Y))E. \nonumber
\end{align*}
So, $D_{0}$ is integrable if and only if the components of $[X,Y]$ with respect to $\bar{\varphi}E,\bar{\varphi}N,E$  vanish, therefore if and only if $C(X,Y)=C(Y,X)$, $ C(X,\bar{\varphi}Y)=C(Y,\bar{\varphi}X)$, $B(X,\bar{\varphi}Y)=B(Y,\bar{\varphi}X)$.
\end{proof}
\begin{remark}\label{simmetria}\emph{
Looking at (\ref{2°forma}) and using the above theorem, we deduce that 
 $\overset{\circ}{h}$ is symmetric on $D_0$ if and only if $D_0$ is integrable. Moreover, the integrability of $D_{0}$ implies that $\overset{\circ}{\nabla}$ is a linear symmetric connection on the integral manifolds.}
\end{remark}
\begin{corollary}\label{C,phi}
If $D_{0}$ is an integrable distribution, then for $X,Y\in\Gamma(D_{0})$ we have
\begin{eqnarray*}
B(\bar{\varphi}X,\bar{\varphi}Y)=-B(X,Y),\quad C(\bar{\varphi}X,\bar{\varphi}Y)=-C(X,Y).
\end{eqnarray*}
\end{corollary}

\medskip
Looking at the decomposition (\ref{decomp06}) and considering the symmetric connection $\nabla$, we can define,
as usual, the unsymmetrized second fundamental
form of $D_{0}$,  $A^{D_{0}}$, setting
\[
A_{X}^{D_{0}}Y=p_{\mathcal{F}}(  \nabla_{X_{0}}Y_{0})  ,
\]
for $X,Y\in\Gamma(  TM)  $, where $X_{0}$, $Y_{0}$ are the
projection of $X$ and $Y$ onto $D_{0}$ and $p_{\mathcal{F}}:TM\rightarrow
\mathcal{F}$ is the canonical projection on $\mathcal{F}$. Then, using (\ref{GaussD0}), since 
$p_{\mathcal{F}}( \overset{\circ}{\nabla}_{X_{0}}Y_{0})=0$, we get 
$A_{X}^{D_{0}}Y=\overset{\circ}{h}(  X_{0},Y_{0})$,
and the symmetrized second fundamental form
$B^{D_{0}}$ is given by
\begin{align*}
B^{D_{0}}(  X,Y) =\frac{1}{2}\{\overset{\circ}{h}(  X_{0},Y_{0})  ) +\overset{\circ}{h}(
Y_{0},X_{0})  )  \},
\end{align*}
for any $X,Y\in\Gamma(  TM)  $. Furthermore the mean curvature vector of the distribution $D_0$, $D_0$ being integrable or not integrable, is defined as
\begin{equation}
\label{media}
\mu^{D_0} = \frac{1}{rank(D_0)}trace(B^{D_0}),
\end{equation}
and $D_0$ is called minimal (respectively totally geodesic) if $\mu^{D_0}$ (respectively $B^{D_0}$) vanishes.  
\begin{proposition}
Let $(  \bar{M},\bar{\varphi},\bar
{\xi}_{\alpha},\bar{\eta}^{\alpha},\bar{g})  $ be an indefinite
$\mathcal{S}$-manifold and $(  M,g,S(  TM)  )  $
a characteristic lightlike hypersurface of $\bar{M}$ such that the distribution $D_{0}$ is integrable. Then, $D_0$ is minimal with respect to the symmetric connection $\nabla$ on $M$ and all its integral manifolds are minimal submanifolds of $M$ with respect to $\nabla$.
\end{proposition}

\begin{proof}
We note that $rank\,(  D_{0})  =2n+r-4=2(  n-2)  +r$,
 and, $D_{0}$ being integrable, $\overset{\circ}{h}$ is symmetric and $\mu^{D_0} = \frac{1}{2(  n-2)  +r}trace(B^{D_0})=\frac{1}{2(  n-2)+r}trace(\overset{\circ}{h}) $.
We consider an adapted frame in $D_0$, 
$( X_{a},\bar{\varphi}X_{a},\bar{\xi}_{\alpha})$ with $ a\in\{1,\ldots,n-2\}$ and $\alpha\in\{1,\ldots,r\}$,
and we have 
\begin{equation*}
trace(  \overset{\circ}{h})  =\sum_{a=1}^{n-2}\varepsilon_{a}(  \overset{\circ}{h}(  X_{a},X_{a})  +\overset{\circ}%
{h}(  \bar{\varphi}X_{a},\bar{\varphi}X_{a})  )
+\sum_{\alpha=1}^{r}\varepsilon_{\alpha}\overset{\circ}{h}(  \bar{\xi
}_{\alpha},\bar{\xi}_{\alpha}).
\end{equation*}
We get $\overset{\circ}{h}(  \bar{\xi}_{\alpha},\bar{\xi}_{\alpha})
=-C(  \bar{\xi}_{\alpha},\bar{\varphi}\bar{\xi}_{\alpha})
\bar{\varphi}E-B(  \bar{\xi}_{\alpha},\bar{\varphi}\bar{\xi}_{\alpha
})  \bar{\varphi}N+C(  \bar{\xi}_{\alpha},\bar{\xi}_{\alpha
})  E=C(  \bar{\xi}_{\alpha},\bar{\xi}_{\alpha})  E=0$,
since (\ref{num05}) and (\ref{equ12}) imply 
$C(\bar{\xi}_{\alpha},\bar{\xi}_{\alpha})= g(A_{N}\bar{\xi}_{\alpha},\bar{\xi}_{\alpha}) = \bar{g}(A_{N}\bar{\xi}_{\alpha},\bar{\xi}_{\alpha})=-\bar{g}(\bar{\nabla}_{\bar{\xi}_{\alpha}}N,\bar{\xi}_{\alpha}) = \bar{g}(N,\bar{\nabla}_{\bar{\xi}_{\alpha}}\bar{\xi}_{\alpha})=0$.
Furthermore, we know that $\bar{\varphi}^{2}X_{a}=-X_{a}$, therefore we
deduce
\begin{align*}
\overset{\circ}{h}(  X_{a},X_{a})  +\overset{\circ}{h}(
\bar{\varphi}X_{a},\bar{\varphi}X_{a}) & =-C(  X_{a},\bar{\varphi}X_{a})  \bar{\varphi}E-B(
X_{a},\bar{\varphi}X_{a})  \bar{\varphi}N +C(  X_{a},X_{a})  E\\
& \quad +C(  \bar{\varphi}X_{a},X_{a})
\bar{\varphi}E +B(  \bar{\varphi}X_{a},X_{a})  \bar{\varphi}N+C(
\bar{\varphi}X_{a},\bar{\varphi}X_{a})  E\\
& =(  C(  X_{a},X_{a})  +C(  \bar{\varphi}X_{a}%
,\bar{\varphi}X_{a})  )  E=0,
\end{align*}
since, using Corollary \ref{C,phi}, we get $C(  X_{a},X_{a})  +C(  \bar{\varphi}X_{a},\bar{\varphi}%
X_{a}) =0$. 
This completes the proof.
\end{proof}

\bigskip
Now, we consider the decomposition (\ref{decomp07})  
and for any $X\in\Gamma(  TM)  $, $Y\in\Gamma(  D_{0})
$ and $W\in\Gamma(  \mathcal{E})  $, we have 
\begin{equation*}
\bar{\nabla}_{X}Y=\widetilde{\nabla}_{X}Y+\widetilde{h}(  X,Y),
\end{equation*}
where $\widetilde{\nabla}$ is a linear
connection on $D_{0}$ and $\widetilde{h}:\,\Gamma(  TM)  \times\Gamma(  D_{0})
\rightarrow\Gamma(  \mathcal{E})$ is $\mathfrak{F}(  M)$-bilinear.

\noindent Let $\mathcal{U}\subset M$ be a coordinate neighbourhood as fixed in Theorem
\ref{teorsezioneN} and $X,Y\in\Gamma(  D_{0}|
_{\mathcal{U}})  $. Then, using (\ref{decomp07}),
 and putting
$F_{1}(  X,Y) =\bar{g}(  \widetilde{h}(  X,Y)
,\bar{\varphi}N)$, $ F_{2}(  X,Y)  =\bar{g}(  \widetilde{h}(  X,Y)
,\bar{\varphi}E)$, $F_{3}(  X,Y) =\bar{g}(  \widetilde{h}(  X,Y)
,N)$, and $F_{4}(  X,Y)  =\bar{g}(  \widetilde{h}(  X,Y)
,E)$,
 we can write locally:
\begin{align}
\bar{\nabla}_{X}Y &=\widetilde{\nabla}_{X}Y+F_{1}(  X,Y)
\bar{\varphi}E+F_{2}(  X,Y)  \bar{\varphi}N+F_{3}(
X,Y)  E+F_{4}(  X,Y)  N. \label{equ39} 	
\end{align}
Now, we express the $F_{i}$'s, $i\in\{  1,2,3,4\}  $, in terms of $B$
and $C$.
\noindent We begin to compute $F_{3}$ and $F_{4}$. For any $X,Y\in\Gamma(
D_{0}|  _{\mathcal{U}})  $, from (\ref{equ39}), using
(\ref{equ12}), (\ref{num05}) and being $\bar{\nabla}$ a metric connection,
we have
\begin{align*}
F_{3}(  X,Y) = -\bar{g}(  Y,\bar{\nabla}_{X}N)  =-\bar{g}(  Y,-A_{N}%
X+\tau(  X)  N) =-\bar{g}(  Y,-A_{N}X)  =g(  A_{N}X,Y)  =C(X,Y),
\end{align*}
and, again from (\ref{equ39}), using (\ref{equ12}), (\ref{equ21}),
(\ref{num05}), we have
\begin{align*}
F_{4}(  X,Y)   & =-\bar{g}(  Y,\bar{\nabla}_{X}E)  =-\bar{g}(  Y,\nabla
_{X}E+B(  X,E)  N)  =-\bar{g}(  Y,\nabla_{X}E)\\  
& =-\bar{g}(  Y,-A_{E}%
X-\tau(  X)  E)  =g(  A_{E}X,Y) =B(  X,Y).
\end{align*}
For $F_{2}$, using Lemma \ref{cond08},
(\ref{equ12}), (\ref{equ21}), (\ref{num05}), and $D_{0}$ being a
$\bar{\varphi} $-invariant distribution, we have:
\begin{align*}
F_{2}(  X,Y)&=\bar{g}(  \bar{\nabla}_{X}Y,\bar{\varphi
}E)  =-\bar{g}(  \bar{\varphi}(  \bar{\nabla}_{X}Y)
,E)=-\bar{g}(  \bar{\nabla}_{X}(  \bar{\varphi}Y)  ,E)
=\bar{g}(  \bar{\varphi}Y,\bar{\nabla}_{X}E)=\bar{g}(  \bar{\varphi}Y,\nabla_{X}E)\\
&=-g(  \overset{*}{A}_{E}X,\bar{\varphi}Y)=-B(X,\bar{\varphi}Y).
\end{align*}
Again, using Lemma \ref{cond08},
(\ref{equ12}) and (\ref{num05}), also by the $\bar{\varphi}$-invariance of
distribution $D_{0}$, we compute:
\begin{align*}
F_{1}(  X,Y)&=\bar{g}(  \bar{\nabla}_{X}Y,\bar{\varphi
}N)  =-\bar{g}(  \bar{\varphi}(  \bar{\nabla}_{X}Y)
,N)=-\bar{g}(  \bar{\nabla}_{X}(  \bar{\varphi}Y)  ,N)
=\bar{g}(  \bar{\varphi}Y,\bar{\nabla}_{X}N)\\
&=-g(  A_{N}X,\bar{\varphi}Y) =-C(  X,\bar{\varphi}Y).
\end{align*}
Then, (\ref{equ39}) becomes
$\bar{\nabla}_{X}Y =\widetilde{\nabla}_{X}Y-C(  X,\bar{\varphi
}Y)  \bar{\varphi}E-B(  X,\bar{\varphi}Y)  \bar{\varphi
}N+C(  X,Y)  E+B(  X,Y)  N$
and, locally, $\widetilde{h}(  X,Y)=-C(  X,\bar{\varphi}Y)
\bar{\varphi}E-B(  X,\bar{\varphi}Y)  \bar{\varphi}N+C(
X,Y)  E+B(  X,Y)  N=\overset{\circ}{h}(X,Y)+B(  X,Y)  N$.
\begin{proposition}
Let $(  \bar{M},\bar{\varphi},\bar{\xi}_{\alpha},\bar{\eta
}^{\alpha},\bar{g})  $ be an indefinite $\mathcal{S}$-manifold and
$(  M,g,S(  TM)  )  $ a characteristic lightlike hypersurface of
$\bar{M}$. Supposing $D_{0}$ integrable, we have
$trace(  \widetilde{h})  =0$,
that is all the integral manifolds of $D_{0}$ are minimal submanifolds of
$\bar{M}$ and $D_{0}$ is minimal.
\end{proposition}

\begin{proof}
We have:
 $trace(  \widetilde{h})  =\sum_{a=1}^{n-2}\varepsilon_{a}(
\widetilde{h}(  X_{a},X_{a})  +\widetilde{h}(  \bar{\varphi
}X_{a},\bar{\varphi}X_{a})  )  +\sum_{\alpha=1}^{r}\varepsilon
_{\alpha}\widetilde{h}(  \bar{\xi}_{\alpha},\bar{\xi}_{\alpha})$.

\noindent We note that $\widetilde{h}(  \bar{\xi}_{\alpha},\bar{\xi}_{\alpha})   =	\overset{\circ}{h}(  \bar{\xi}_{\alpha},\bar{\xi}_{\alpha})+B(  \bar{\xi}_{\alpha},\bar{\xi}_{\alpha})  N= \overset{\circ}{h}(  \bar{\xi}_{\alpha},\bar{\xi}_{\alpha})$,
since, using (\ref{num05}), (\ref{equ21}) and (\ref{equ12}),  
\begin{align*}
B(  \bar{\xi}_{\alpha},\bar{\xi}_{\alpha})   & =g(
\overset{*}{A}_{E}\bar{\xi}_{\alpha},\bar{\xi}_{\alpha})  =\bar
{g}(  \overset{*}{A}_{E}\bar{\xi}_{\alpha},\bar{\xi}_{\alpha}) =-\bar{g}(  \nabla_{\bar{\xi}_{\alpha}}E+\tau(  \bar{\xi}_{\alpha
})  E,\bar{\xi}_{\alpha}) =-\bar{g}(  \nabla_{\bar{\xi}_{\alpha}}E,\bar{\xi}_{\alpha}) \\
& =-\bar{g}(  \bar{\nabla}_{\bar{\xi}_{\alpha}}E-B(  \bar{\xi
}_{\alpha},E)  N,\bar{\xi}_{\alpha}) =-\bar{g}(  \bar{\nabla}_{\bar{\xi}_{\alpha}}E,\bar{\xi}_{\alpha
})  =\bar{g}(  E,\bar{\nabla}_{\bar{\xi}_{\alpha}}\bar{\xi}%
_{\alpha})  =0.
\end{align*}
Being $\widetilde{h}(  X_{a},X_{a})  +\widetilde{h}(  \bar{\varphi
}X_{a},\bar{\varphi}X_{a}) =	\overset{\circ}{h}(  X_{a},X_{a})+B(  X_{a},X_{a})  N
+ \overset{\circ}{h}(  \bar{\varphi
}X_{a},\bar{\varphi}X_{a})+B(  \bar{\varphi
}X_{a},\bar{\varphi}X_{a})  N$, using Corollary \ref{C,phi}, we obtain
$trace(  \widetilde{h})  =trace(\overset{\circ}{h})=0$.
\end{proof}
\begin{proposition}
Let $(  \bar{M},\bar{\varphi},\bar{\xi}_{\alpha},\bar{\eta
}^{\alpha},\bar{g})  $ be a indefinite $\mathcal{S}$-manifold of dimension $2n+r$, $n\geq 3$ and
$(  M,g,S(  TM)  )  $ a characteristic lightlike
hypersurface of $\bar{M}$. If $D_{0}$ is an integrable distribution, then the
leaves of $D_{0}$ have an indefinite $\mathcal{S}$-structure.
\end{proposition}%
\begin{proof}
Let $M_{0}$ be a leaf of $D_{0}$, then for any $p\in M_{0}$ we have
$T_{p}M_{0}=(  D_{0})  _{p}$ and $dim M_0=2(n-2)$. If $X_{0}\in TM_{0}$, we have
\[
\varphi X_{0}=\bar{\varphi}SX_{0}=\bar{\varphi}X_{0},
\]
being $S:\Gamma(  TM)  \rightarrow\Gamma(  D)  $
and $D=D_{0}\bot\bar{\varphi}(  Rad(TM))  \bot Rad(TM)$. We put $\overset{\circ}{\varphi:}=\varphi|_{D_{0}}$, and for any $\alpha\in\{  1,\ldots,r\}  $ $\overset{\circ}{\eta}^{\alpha}:=\bar{\eta}^{\alpha}|  _{D_{0}}$ so 
$\overset{\circ}{\varphi}$ defines an $(  1,1)  $-type tensor field
on $M_{0}$ because $D_{0}$ is $\bar{\varphi}$-invariant. Now we consider
$(  M_{0},\overset{\circ}{\varphi},\bar{\xi}_{\alpha},\bar{\eta}^{\alpha
},g)  $ and check that this is an indefinite $\mathcal{S}$-structure.
We know that $\varphi^{2}X=-X+\bar{\eta
}^{\alpha}(  X)  \bar{\xi}_{\alpha}+u(  X)  U$, for any $X\in\Gamma( TM)$, and that
$u(  Y)  =0$ for any $Y\in\Gamma( D)$, so we deduce
\[
\overset{\circ}{\varphi}^{2}X_{0}=-X_{0}+\overset{\circ}{\eta}^{\alpha}(
X_{0})  \bar{\xi}_{\alpha},
\]
for any $X_{0}\in \Gamma(TM_{0})$.
For any $\alpha,\beta\in\{  1,\ldots,r\}  $ we have $\overset{\circ}{\eta}^{\alpha}(  \bar{\xi}_{\beta})  =\delta
_{\beta}^{\alpha}$ and then $(  M_{0},\overset{\circ}{\varphi},\bar{\xi}_{\alpha},\overset
{\circ}{\eta}^{\alpha})  $ is a $g.f.f$-manifold. Now, to prove
the compatibility between the $g.f.f$-structure and the metric $g$ on $M_{0}$,
by the definition in (\ref{equ27}), for any
$X_{0},Y_{0}\in\Gamma( TM_{0})$ we have
\begin{align*}
g(  \overset{\circ}{\varphi}X_{0},\overset{\circ}{\varphi}Y_{0})
& =\bar{g}(  \bar{\varphi}X_{0},\bar{\varphi}Y_{0})  -u(
Y_{0})  \bar{g}(  \bar{\varphi}X_{0},N) -u(  X_{0})  \bar{g}(  N,\bar{\varphi}Y_{0})
+u(  X_{0})  u(  Y_{0})  \bar{g}(  N,N)\\
&  =g(  X_{0},Y_{0})  -\varepsilon_{\alpha
}\overset{\circ}{\eta}^{\alpha}(  X_{0})  \overset{\circ}{\eta}^{\alpha}(  Y_{0}).
\end{align*}
Moreover, for any $X_{0},Y_{0}\in\Gamma(TM_{0})$, for any $\alpha\in\{1,\ldots,r\}$, we get
$d\bar{\eta}^{\alpha}(X_{0},Y_{0})=d\overset{\circ}{\eta}^{\alpha}( X_{0}, Y_{0})$, 
and $\overset{\circ}{\Phi}(X_{0},Y_{0})=g(X_{0},\overset{\circ}{\varphi}Y_{0})=\bar{g}(X_{0},\bar{\varphi}Y_{0})=d\bar{\eta}^{\alpha}(X_{0},Y_{0})=d\overset{\circ}{\eta
}^{\alpha}( X_{0}, Y_{0})$. Finally, $(M_{0},\overset{\circ}{\varphi},\bar{\xi}_{\alpha},\overset
{\circ}{\eta}^{\alpha},g)$ is an indefinite $\mathcal{S}$-manifold, since
 $N_{\overset{\circ}{\varphi}}+2d\overset{\circ}{\eta
}^{\alpha}\otimes \bar\xi_{\alpha}$ and $ N_{\bar{\varphi}}+2d\bar\eta^{\alpha}\otimes \bar\xi_{\alpha}$ coincide on $D_0$.
Moreover, $\overset{\circ}{\nabla}$ is the Levi-Civita connection on $M_{0}$. In fact, by Remark \ref{remNab}, since $D_0\subset S(TM)$, for any $X_0,Y_0,Z_0\in\Gamma{T M_{0}}$ we have
\begin{align*}
(  \overset{\circ}{\nabla}_{X_{0}}g)(
Y_{0},Z_0)&=X_0(g(Y_0,Z_0))-g( \overset{\circ}{\nabla}_{X_{0}} Y_{0},Z_{0})-g( Y_{0}, \overset{\circ}{\nabla}_{X_{0}}Z_{0})=X_0(g(Y_0,Z_0))-g({\nabla}_{X_{0}} Y_{0},Z_{0})\\
&\quad +g( \overset{\circ}{h}(X_{0}, Y_{0}),Z_{0})-g( Y_{0},{\nabla}_{X_{0}}Z_{0})+g( Y_{0}, \overset{\circ}{h}(X_{0},Z_{0}))=({\nabla}_{X_{0}}g)(
Y_{0},Z_0)=0.
\end{align*}
Hence $\overset{\circ}{\nabla}$ is a metric connection. By Remark \ref{simmetria} it is also symmetric, thus it is the Levi-Civita connection and, from (\ref{LC}), we have $(\overset{\circ}{\nabla}_{X_0}\overset{\circ}{\varphi})Y_0=g(\overset{\circ}{\varphi}X_0,\overset{\circ}{\varphi}Y_0)\bar{\xi}+\bar{\eta}(Y_0)\overset{\circ}{\varphi}^{2}(X_0)$.
 \end{proof}
\subsection{The distribution $D$}
By a direct computation, one can prove the following result.
\begin{lemma}
\label{lemma agg}Let $(\bar{M},\bar{\varphi},\bar{\xi}_{\alpha},\bar{\eta}^{\alpha},\bar{g})$ be an indefinite $\mathcal{S}$-manifold and consider $(M,g,S(TM))$ a characteristic lightlike hypersurface. Then the component of $(\bar{\nabla}_{X}\bar{\varphi})Y$ along $ltr(TM)$ vanishes, for any $X\in\Gamma(TM)$ and $Y\in\Gamma(T\bar{M})$.
\end{lemma}
\begin{proposition}
Let $(\bar{M},\bar{\varphi},\bar{\xi}_{\alpha},\bar{\eta}^{\alpha},\bar{g})$ be an indefinite $\mathcal{S}$-manifold and $(M,g,S(TM))$ a characteristic lightlike hypersurface. Then $D=D_{0}\bot\bar{\varphi}(Rad(TM))\bot Rad(TM)$ is integrable if and only if $B$ satisfies the following conditions:
\begin{itemize}
    \item[a)]$B(X,\bar{\varphi}Y)=B(\bar{\varphi}X,Y)$, for any $X,Y\in\Gamma(D_{0})$
    \item[b)]$B(X,V)=0$, for any $X\in\Gamma(D_{0})$
    \item[c)]$B(V,V)=0$,
\end{itemize}
where $V=-\bar{\varphi}E$.
\end{proposition}
\begin{proof}
At first, for any $X,Y\in\Gamma(D)$ we compute the component of $[X,Y]$ along $\bar{\varphi}(ltr(TM))$
\begin{align*}
    \bar{g}([X,Y],\bar{\varphi}E)& = -\bar{g}(\bar{\varphi}\bar{\nabla}_{X}Y,E)+\bar{g}(\bar{\varphi}\bar{\nabla}_{Y}X,E) = -\bar{g}(\bar{\nabla}_{X}(\bar{\varphi}Y),E) +\bar{g}(\bar{\nabla}_{Y}(\bar{\varphi}X),E)
    \\
    & = \bar{g}(\bar{\varphi}Y,\bar{\nabla}_{X}E)-\bar{g}(\bar{\varphi}X,\bar{\nabla}_{Y}E)= -g(\bar{\varphi}Y,\overset{*}{A}_{E}X)+g(\bar{\varphi}X,\overset{*}{A}_{E}Y).
\end{align*}
From the definition of $D$ we get $X = \alpha E+ \beta \bar{\varphi}E+ X_{0}$, and $ Y = \delta E+ \gamma \bar{\varphi}E+ Y_{0}$.
Using the previous expression of $X,Y\in\Gamma(D)$, being $D$ $\bar{\varphi}$-invariant and $B(E,X)=0$ for any $X\in\Gamma(TM)$, we have
\begin{align}
    \bar{g}([X,Y],\bar{\varphi}E)& = (\gamma\alpha - \beta\delta) B(V,V)- \gamma B(V,\bar{\varphi}X_{0}) - \alpha B(Y_{0},V)\label{eqrif}\\
    & \quad +\beta B(V,\bar{\varphi}Y_{0})+\delta B(X_{0},V) + B(Y_{0},\bar{\varphi}X_{0}) -B(\bar{\varphi}Y_{0},X_{0})\nonumber.
\end{align}
So, if we suppose that $D$ is integrable, being $\bar{\varphi}E$, $E$, $X_{0}$ and $Y_{0}$ sections of $D$, then we get $0=\bar{g}([\bar{\varphi}E,E],\bar{\varphi}E) = -\bar{g}(\bar{\varphi}E,\bar{\varphi}E)=-B(V,V)$.
Finally, if $X\in\Gamma(D_{0})$ we find $0 = \bar{g}([X,E],\bar{\varphi}E) = B(E,\bar{\varphi}X)-B(X,\bar{\varphi}E)=B(X,V)$, and if $X,Y\in\Gamma(D_{0})$ we get $0 = \bar{g}([X,Y],\bar{\varphi}E)=B(\bar{\varphi}Y,X)-B(Y,\bar{\varphi}X)$.
Vice versa, using (\ref{eqrif}) and a), b), c), it is easy to check that $[X,Y]$ in $\Gamma(D)$.
\end{proof}

Now, a consequence is the following proposition.
\begin{proposition}
Let $(\bar{M},\bar{\varphi},\bar{\xi}_{\alpha},\bar{\eta}_{\alpha},\bar{g})$ be an indefinite $\mathcal{S}$-manifold and $(M,g,S(TM))$ a characteristic lightlike hypersurface. If $(M,g,S(TM))$ is totally geodesic, then the following statements hold:
\begin{itemize}
    \item [a)] the distribution $D$ is integrable;
    \item [b)] the distribution $D$ is parallel with respect to the induced connection $\nabla$;
    \item [c)] $M$ is locally a product $M^{*}\times C$, where $M^{*}$ is a leaf of $D$ and $C$ is a lightlike curve tangent to the distribution $\bar{\varphi}(ltr(TM))$.
\end{itemize}
\end{proposition}
\begin{proof}
Being $(M,g,S(TM))$ totally geodesic, a) follows from the previous proposition.

\noindent To prove b), we need only to check $g(\nabla_{X}E,\bar{\varphi}E)=0$, $g(\nabla_{X}\bar{\varphi}E,\bar{\varphi}E)=0$ and $g(\nabla_{X}Y_{0},\bar{\varphi}E)=0$ for any $X\in\Gamma(TM)$ and $Y_{0}\in\Gamma(D_{0})$. Hence, using Lemma \ref{lemma agg}, we get
\begin{align*}  
g(\nabla_{X}E,\bar{\varphi}E)&=\bar{g}(\bar{\nabla}_{X}E,\bar{\varphi}E)=-\bar{g}(E,\bar{\nabla}_{X}\bar{\varphi}E)=-B(X,\bar{\varphi}E)=0,\\
g(\nabla_{X}\bar{\varphi}E,\bar{\varphi}E) & =\bar{g}(\bar{\nabla}_{X}\bar{\varphi}E,\bar{\varphi}E)=-\bar{g}(\bar{\varphi}\bar{\nabla}_{X}\bar{\varphi}E,E)=\bar{g}((\bar{\nabla}_{X}\varphi)\bar{\varphi}E,E)\\
& \quad -\bar{g}(\bar{\nabla}_{X}\bar{\varphi}^{2}E,E)=\bar{g}(\nabla_{X}E,E)=0,\\
g(\nabla_{X}Y_{0},\bar{\varphi}E) & =\bar{g}(\bar{\nabla}_{X}Y_{0},\bar{\varphi}E)=-\bar{g}(\bar{\varphi}\bar{\nabla}_{X}Y_{0},E)=\bar{g}((\bar{\nabla}_{X}\bar{\varphi})Y_{0},E)\\
& \quad -\bar{g}(\bar{\nabla}_{X}(\bar{\varphi}Y_{0}),E)= -\bar{g}(\bar{\nabla}_{X}(\bar{\varphi}Y_{0}),E)= -B(X,\bar{\varphi}Y_{0})=0.
\end{align*}
Finally, from a) we deduce that $D$ determines a foliation. Being $\bar{\varphi}(ltr(TM))$ a 1-dimensional distribution, it defines a foliation. So, being $TM=D\oplus\bar{\varphi}(ltr(M))$, we obtain c).
\end{proof}
\section{Totally umbilical lightlike hypersurface and totally umbilical screen distribution}

We begin to prove the following Lemma.
\begin{lemma}\label{R e R bar}
Let  $( \bar{M},\bar{\varphi
},\bar{\xi}_{\alpha},\bar{\eta}^{\alpha},\bar{g})  $ be an indefinite $\mathcal{S}$-manifold and $(M,g,S(TM))$ a lightlike hypersurface. Then, the Riemannian $(0,4)$-type curvature tensor fields $\bar R$ and $R$ of $\bar{M}$ and $M$ verify the following relations, for any $X,Y,Z\in \Gamma(TM)$,
\begin{align*}
\bar{R}(X,Y,Z,E)&=-\{(\nabla_{X}B)(Y,Z)-(\nabla_{Y}B)(X,Z)+\tau(X)B(Y,Z)-\tau(Y)B(X,Z)\},\\
\bar{R}(X,Y,Z,N)&=R(X,Y,Z,N),
\end{align*}
 \end{lemma}

\begin{proof}
Only using the Gauss and Weingarten equations for lightlike hypersurfaces, we get
\begin{align*}
	\bar{R}(X,Y,Z,E)& =-\{B(X,\nabla_{Y}Z)+X(B(Y,Z))+\tau(X)B(Y,Z)-B(Y,\nabla_{X}Z)\\
	& \quad -Y(B(X,Z))-\tau(Y)B(X,Z)+B(\nabla_{Y}X,Z)-B(\nabla_{X}Y,Z)\}\\
	&=-\{(\nabla_{X}B)(Y,Z)-(\nabla_{Y}B)(X,Z)+\tau(X)B(Y,Z) -\tau(Y)B(X,Z)\},\\
	\bar{R}(X,Y,Z,N)&=-\bar{g}(\bar{R}(X,Y,Z),N)=-\{\bar{g}(\bar{\nabla}_{X}\nabla_{Y}Z,N)
	+B(Y,Z)\bar{g}(\bar{\nabla}_{X}N,N)\\
&\quad -\bar{g}(\bar{\nabla}_{Y}\nabla_{X}Z,N)+B(X,Z)\bar{g}(\bar{\nabla}_{Y}N,N)-\bar{g}(\nabla_{[X,Y]}Z,N)\}\\
	&=-\{\bar{g}(\nabla_{X}\nabla_{Y}Z,N) -\bar{g}(\nabla_{Y}\nabla_{X}Z,N)-\bar{g}(\nabla_{[X,Y]}Z,N)\}\\
	&=R(X,Y,Z,N).
\end{align*}
\end{proof}
\begin{definition}\emph{
Let $(\bar{M},\bar{\varphi},\bar{\xi}_{\alpha},\bar{\eta}^{\alpha},\bar{g})$ be an indefinite $\cal{S}$-manifold and $(M,g,S(TM))$ a lightlike hypersurface. Then $M$ is called \emph{totally umbilical} if, for any coordinate neighbourhood $\cal{U}$, there exists a function $\rho\in\mathfrak{F}(\mathcal{U})$ such that $B(X,Y)=\rho g(X,Y)$ for any $X,Y\in\Gamma{(TM_{|\cal{U}})}$.}
\end{definition}
\begin{theorem}
Let $(\bar{M}(c),\bar{\varphi},\bar{\xi}_{\alpha},\bar{\eta}^{\alpha},\bar{g})$ be an indefinite $\cal{S}$-space form and $(M,g,S(TM))$ a characteristic lightlike hypersurface. If $(M,g,S(TM))$ is totally umbilical then $c=\varepsilon=\sum_{\alpha=1}^{r}\varepsilon_{\alpha}$.
\end{theorem}

\begin{proof}
Being $\bar{M}(c)$ an indefinite $\mathcal{S}$-space form, the Riemannian curvature ${\bar R}$ is given by
\begin{align*}
	{\bar R}(  X,Y,Z,W)   &  =-\frac{c+3\varepsilon}{4}\{  
{\bar g}({\bar\varphi}Y,{\bar\varphi}Z)  {\bar g}({\bar\varphi} X,{\bar\varphi}W) -{\bar g}(  {\bar\varphi}X,{\bar\varphi}Z)
{\bar g}( {\bar\varphi}Y,{\bar\varphi}W)\} \\
&\quad-\frac{c-\varepsilon}{4}\{  \Phi(  W,X)
\Phi(  Z,Y)     -\Phi(  Z,X)  \Phi(  W,Y) +2\Phi(  X,Y)  \Phi(  W,Z)
\} \\
&  \quad -\{ \bar {\eta}(  W)  \bar {\eta
}(  X)  {\bar g}(  {\bar\varphi}Z,{\bar\varphi
}Y)  -\bar {\eta}(  W)  \bar {\eta}(
Y)  {\bar g}(  {\bar\varphi}Z,{\bar\varphi}X)  +\bar {\eta
}(  Y)  \bar {\eta}(  Z)  {\bar g}(
{\bar\varphi}W,{\bar\varphi}X)\\
&  \quad    -\bar {\eta}(  Z)  \bar {\eta}(  X)  {\bar g}(  {\bar\varphi}W,{\bar\varphi}Y)
\},
\end{align*}
 for any $X,Y,Z,W \in \Gamma(T\bar{M})$, (\cite{LP}).
 Since $\bar{g}(\bar{\varphi}E,\bar{\varphi}X)=0$, for any $X\in\Gamma(TM)$, and $\bar\eta (E)=0$, then for any $X,Y,Z \in\Gamma(TM)$,  we get
\begin{align*}
\bar{R}(  X,Y,Z,E)   &  =-\frac{c-\varepsilon}{4}\{  \Phi(  E,X)
\Phi(  Z,Y) -\Phi(  Z,X)  \Phi(  E,Y) +2\Phi(  X,Y)  \Phi(  E,Z)\}\\
&=-\frac{c-\varepsilon}{4}\{  \bar{g}(V,X)\Phi(  Z,Y) -\Phi(  Z,X)  \bar{g}( V,Y) +2\Phi(  X,Y)  \bar{g}(  V,Z)\}. \nonumber 
\end{align*}
On the other hand, from Lemma \ref{R e R bar} we deduce
\begin{align*}
\bar{R}(X,Y,Z,E)=-\{(\nabla_{X}B)(Y,Z)-(\nabla_{Y}B)(X,Z)+\tau(X)B(Y,Z)-\tau(Y)B(X,Z)\}.
\end{align*}
So, replacing $X,Y,Z$ by $PX,E,PZ$ in the above two equations, $P$ being the projection on $S(TM)$ with respect to the decomposition (\ref{decomp01}), we find
\begin{align*}
\bar{R}(  PX,E,PZ,E) =-\frac{c-\varepsilon}{4}\{  -\bar{g}(V,PX)\bar{g}(  PZ,V)-2\bar{g}(  X,V)  \bar{g}(  V,Z)\}=\frac{3}{4}(c-\varepsilon)u(PZ)u(PX).
\end{align*}
On the other hand we have 
\begin{align*}
\bar{R}(PX,E,PZ,E)&=-\{-B(\nabla_{PX}E,PZ)-E(B(PX,PZ))+B(\nabla_{E}PX,PZ)+B(PX,\nabla_{E}PZ)\\
& \quad -\tau(E)B(PX,PZ)\}.
\end{align*}
Then, comparing the two equations and using Gauss and Weingarten equations, we get
\begin{align*}
\frac{3}{4}(c-\varepsilon)u(PZ)u(PX)&=-\{\rho B(PX,PZ)-E(\rho)g(PX,PZ)-\rho(B(E,PX)\bar{g}(PZ,N)\\
&\quad -B(E,PZ)\bar{g}(PX,N))-\rho\tau(E)g(PX,PZ)\}.
\end{align*}
Therefore, being $B(X,E)=0$ for any $X\in\Gamma(TM)$, we still find 
\[
    \frac{3}{4} (c-\varepsilon)u(PZ)u(PX) = -(\rho^{2}-E(\rho)-\rho\tau(E)) \bar{g}(PX,PZ).
\]
Choosing $X=Z=U\in S(TM)$, we have $PX=PZ=U$ and being $u(U)=1$ and $g(U,U)=0$, from the above equation we obtain
$c=\varepsilon$.
\end{proof}
\begin{corollary}
Let $(\bar{M}(c),\bar{\varphi},\bar{\xi}_{\alpha},\bar{\eta}_{\alpha},\bar{g})$ be an indefinite $\mathcal{S}$-space form. If $c\neq\varepsilon$, then there exists no totally umbilical characteristic lightlike hypersurface.
\end{corollary}
\begin{definition}\emph{
Let $(\bar{M},\bar{\varphi},\bar{\xi}_{\alpha},\bar{\eta}_{\alpha},\bar{g})$ be an indefinite $\cal{S}$-manifold and $(M,g,S(TM))$ a lightlike hypersurface. The screen distribution $S(TM)$ is called \emph{totally umbilical} if for any coordinate neighbourhood $\cal{U}$ there exists a function $\lambda\in\mathfrak{F}(\cal{U})$ such that $C(X,PY)=\lambda g(X,PY)$, for any $X,Y\in\Gamma{(TM_{|\cal{U}})}$.}
 \end{definition}
\begin{proposition}
Let $(M,g,S(TM))$ be a characteristic lightlike hypersurface of an indefinite $\cal{S}$-space form. If $S(TM)$ is totally umbilical, then $S(TM)$ is totally geodesic.
\end{proposition}

\begin{proof}
From Lemma \ref{R e R bar}, for any $X,Y,Z\in\Gamma(TM)$ we get $\bar{g}(\bar{R}(X,Y,Z),N)=\bar{g}(R(X,Y,Z),N)$.
Then, replacing $Z$ by $PZ$, we have
\begin{align*}
\bar{g}(\bar{R}(X,Y,PZ),N)&=\bar{g}(\nabla_{X}\nabla_{Y}PZ,N)-\bar{g}(\nabla_{Y}\nabla_{X}PZ,N)-\bar{g}(\nabla_{[X,Y]}PZ,N)\\
&=\bar{g}(\overset{*}{\nabla}_{X}\overset{*}{\nabla}_{Y}PZ,N)+C(X,\overset{*}{\nabla}_{Y}PZ)+X(C(Y,PZ))
 -C(Y,PZ)\bar{g}(\overset{*}{A}_{E}X,N)\\
&\quad -C(Y,PZ)\tau(X)-\bar{g}(\overset{*}{\nabla}_{Y}\overset{*}{\nabla}_{X}PZ,N)-C(Y,\overset{*}{\nabla}_{Y}PZ)-Y(C(X,PZ)) \\
&\quad +C(X,PZ)\bar{g}(\overset{*}{A}_{E}Y,N)+C(X,PZ)\tau(Y) +C(\nabla_{Y}X,PZ)-C(\nabla_{X}Y,PZ)\\
&=X(C(Y,PZ))-C(Y,\overset{*}{\nabla}_{X}PZ)-C(\nabla_{X}Y,PZ) -\tau(X)C(Y,PZ)\\
&\quad -Y(C(X,PZ))+C(X,\overset{*}{\nabla}_{Y}PZ)+C(\nabla_{Y}X,PZ)+\tau(Y)C(X,PZ).
\end{align*}
Now, replacing $X$ by $E$ and $Y,Z$ by $U$ in the previous equation, we obtain
\begin{align*}    
\bar{g}(\bar{R}(E,U,U),N)&=E(C(U,U))-C(U,\overset{*}{\nabla}_{E}U)-C(\nabla_{E}U,U)-\tau(E)C(U,U)\\
&\quad -U(C(E,U)) +C(E,\overset{*}{\nabla}_{U}U) +C(\nabla_{U}E,U)+\tau(U)C(E,U)\\
&=-\lambda (g(U,\overset{*}{\nabla}_{E}U)+g(\nabla_{E}U,U)+g(\nabla_{U}E,U)).
\end{align*}
Note that $\overset{*}{\nabla}_{E}U=\nabla_{E}U-C(E,U)E=\nabla_{E}U$ and $\overset{*}{\nabla}_{U}U=\nabla_{U}U-C(U,U)E=\nabla_{U}U$. So, we find
\begin{align*}
    \bar{g}(\bar{R}(E,U,U),N)&=-\lambda(g(U,\nabla_{U}E)+2g(\nabla_{E}U,U))-\lambda \bar{g}(\bar{\nabla}_{U}E,U)=-\lambda \bar{g}(\bar{\varphi}\bar{\nabla}_{U}E,N)\\
    &=-\lambda\bar{g}(\bar{\nabla}_{U}\bar{\varphi}E,N)=\lambda\bar{g}(\bar{\varphi}E,\bar{\nabla}_{U}N)=-\lambda g(\bar{\varphi}E,A_{N}U)\\
&=-\lambda C(U,\bar{\varphi}E)=\lambda^2 g(U,V)=\lambda^2.
\end{align*}
Being $\bar{M}(c)$ an indefinite $\mathcal{S}$-space form, for the Riemannian curvature tensor field we have
\begin{align*}
    \bar{g}(\bar{R}(E,U,U),N) =\frac{c-\varepsilon}{4}\{  -\Phi(  N,E)
 \bar{g}(  \bar{\varphi}N,N)-\Phi(  U,E)   \bar{g}(   U,U)  +2\Phi(  E,U)   \bar{g}( U,U)\}=0 
\end{align*}
So, we obtain $\lambda=0$ and $S(TM)$ is totally geodesic.
\end{proof}
\section{An example of a characteristic lightlike hypersurface}\label{example}
We begin with a general remark.
Consider a hypersurface $M$ of a semi-Riemannian manifold $(\bar{M}^{m},\bar{g})$ locally given by
\[
x^{A}=f^{A}(u^{1},\ldots,u^{m-1}), \qquad rank(\frac{\partial f^{A}}{\partial u^{a}})=m-1,\quad
A\in\{1,\ldots,m\},\quad a\in\{1,\ldots,m-1\}.\]
A natural basis of the tangent bundle of $M$ is
$\frac{\partial}{\partial u^{a}}=\frac{\partial f^{A}}{\partial u^{a}}\frac{\partial}{\partial x^{A}}$ 
and the induced metric on $M$, denoted by $g$, has components given by
$g_{ab}=\bar{g}_{AB}\frac{\partial f^{A}}{\partial u^{a}}\frac{\partial f^{B}}{\partial u^{b}}$.
Obviously, a hypersurface $M$ of $(\bar{M}^{m},\bar{g})$ is lightlike if and only if $rank(g_{ab})\leq m-2$, i.e. 
\begin{equation}
\Delta=det(\frac{\partial f^{A}}{\partial u^{a}}\bar{g}_{AB}\frac{\partial f^{B}}{\partial u^{b}})=0\label{prop1}.
\end{equation}
\begin{proposition}[\cite{K}] \label{prop2}
Let $A\!=\!(a_{ij})\!\in\! M_{m,n}(\mathbb{R})$ and $B\!=\!(b_{jk})\!\in \!M_{n,m}(\mathbb{R})$. If we consider $C=AB$, then we have
\begin{equation*}
det C=\sum_{1\leq j_{1}\leq j_{2}\ldots\leq j_{n}\leq m} \left|
\begin{array}{cccc}
 a_{1j_{1}} & a_{2j_{1}} & \ldots & a_{nj_{1}}\\
 a_{1j_{2}} & a_{2j_{2}} & \ldots & a_{nj_{2}}\\
 \ldots     & \ldots     & \ldots & \ldots    \\
 a_{1j_{n}} & a_{2j_{n}} & \ldots & a_{nj_{n}} 
\end{array}
\right|
\left|
\begin{array}{cccc}
 b_{j_{1}1} & b_{j_{1}2} & \ldots & b_{j_{1}1}\\
 b_{j_{1}1} & b_{j_{1}1} & \ldots & b_{j_{1}1}\\
 \ldots     & \ldots     & \ldots & \ldots    \\
 b_{j_{1}1} & b_{j_{1}1} & \ldots & b_{j_{1}1} 
\end{array}
\right|.
\end{equation*}
\end{proposition}
\bigskip
Now, we refer to the three examples of indefinite $S$-manifolds given in \cite{LP} to look for a lightlike hypersurface. We begin with Example 4.2 in \cite{LP}, $\bar M =(\mathbb{R}^{6}_{2},\varphi,\xi_{\alpha},\eta^{\alpha},g)$, where $(x^1,x^2,y^1,y^2,z^1,z^2)$ are standard coordinates in ${\mathbb R} ^6$,  
$\xi_{\alpha}:=\frac{\partial}{\partial z^{\alpha}}$, $\eta^{\alpha}:=dz^{\alpha}-\sum_{i=1}^{2} \tau_{i}y^{i}dx^{i}$, $\alpha \in \{1,2\}$,
and $\varphi$, $g$ are given by
\begin{eqnarray*}
\begin{array}{ll}
F &=\left(
\begin{array}{ccc}
0    & I_{2}     & 0 \\
-I_{2} & 0  & 0\\
0  & Y  & 0
\end{array}
\right),\qquad {\rm where} \quad
Y=
\left(
\begin{array}{cc}
-y^{1} & y^{2} \\
 -y^{1} & y^{2}
\end{array}
\right),\\
\\
g&=\sum\nolimits_{\alpha=1}^{2}\eta^{\alpha}\otimes\eta^{\alpha}+\frac{1}{2}\sum\nolimits_{i=1}^{2}\tau_{i}((dx^{i})^{2}+(dy^{i})^{2}),
\end{array}
\end{eqnarray*}
respectively, being $\tau_{i}=\mp1$ according to whether $i=1$ or $i=2$. Note that the Killing characteristic vector fields are both spacelike.

We compute all the non null determinants $M_{AB}$, obtained by taking the line $A$ and the column $B$ out of $G$, matrix of $g$, 
\begin{align*}
M_{11}&=M_{33}=-\frac{1}{8},\quad M_{22}=M_{44}=\frac{1}{8}, \quad 	    M_{55}=M_{66}=\frac{1}{16}+\frac{1}{8}(y^{2})^{2}-\frac{1}{8}(y^{1})^{2},\\	
M_{15}&=\frac{1}{8}y^{1},\quad M_{16}=-\frac{1}{8}y^{1}, \quad M_{25}=-\frac{1}{8}y^{2},\quad M_{26}=\frac{1}{8}y^{2},
\quad M_{56}=-\frac{1}{8}(y^{2})^{2}+\frac{1}{8}(y^{1})^{2}.
\end{align*}
Now, in terms of $D^{A}$ and $M_{AB}$, where $D=(\frac{\partial f^{A}}{\partial u^{i}})$ and, for any $A\in\{1,\ldots,6\}$, $D^{A}$ denotes the determinant of the matrix obtained by $D$ deleting the column of index $A$, we state the following result.
\begin{proposition}\label{det1}
Let $M$ be a hypersurface of $\bar M$. Then, $M$ is lightlike if and only if
\begin{align*}
\frac{1}{2}(&-(D^{1})^{2}+(D^{2})^{2}-(D^{3})^{2}+(D^{4})^{2}+(\frac{1}{2}+(y^{2})^{2}-(y^{1})^{2})((D^{5})^{2}+(D^{6})^{2}))\\
& +y^{1}(D^{1}D^{5}-D^{1}D^{6})-y^{2}(D^{2}D^{5}-D^{2}D^{6})+((y^{1})^{2}-(y^{2})^{2})D^{5}D^{6}=0.
\end{align*}
\end{proposition}
\begin{proof}
From (\ref{prop1}) we know that $M$ is a lightlike hypersurface if and only if $\Delta=det(\frac{\partial f^{A}}{\partial u^{a}}g_{AB}\frac{\partial f^{B}}{\partial u^{b}})=0$, thus, using Proposition \ref{prop2}, we have
\begin{align*}
	\Delta &=-\frac{1}{8}((D^{1})^{2}+(D^{3})^{2})+\frac{1}{4}y^{1}(D^{1}D^{5}-D^{1}D^{6})\\
	& \quad +\frac{1}{8}((D^{2})^{2}+(D^{4})^{2})-\frac{1}{4}y^{2}(D^{2}D^{5}-D^{2}D^{6})+\frac{1}{4}((y^{1})^{2}-(y^{2})^{2})D^{5}D^{6}\\
	& \quad +\frac{1}{8}(\frac{1}{2}+(y^{2})^{2}-(y^{1})^{2})((D^{5})^{2}+(D^{6})^{2}).
\end{align*}
This ends the proof.
\end{proof}

Now, we give a condition which  ensures that the characteristic vector fields belong to $TM$.
\begin{proposition}\label{det2}
Let $M$ be a hypersurface of $\bar M$. Then, the characteristic vector fields are tangent to $M$ if and only if $\;D^{5}=D^{6}=0$.
\end{proposition}
\begin{proof}
We suppose that, for any $\alpha\in\{1,2\}$, $\xi_{\alpha}$ is tangent to $M$, then we can write
$\xi_{\alpha}=X_{\alpha}^{a}\frac{\partial}{\partial u^{a}}$.
Being $\frac{\partial}{\partial u^{a}}=\frac{\partial f^{A}}{\partial u^{a}}\frac{\partial}{\partial x^{A}}$, we have
$\xi_{\alpha}=X_{\alpha}^{a}\frac{\partial f^{A}}{\partial u^{a}}\frac{\partial}{\partial x^{A}}$ and,
on the other hand, we know
$\xi_{\alpha}=\frac{\partial}{\partial z^{\alpha}}$.
So, we obtain
\begin{equation*}
{\rm I)}: \left\{
\begin{array}{cc}
	X_{1}^{a}\frac{\partial f^{A}}{\partial u^{a}}=0& \qquad A\neq5\\
	\\
	X_{1}^{a}\frac{\partial f^{5}}{\partial u^{a}}=1&
\end{array}
\right. \qquad {\rm and}\qquad
 {\rm II)}: \left\{
\begin{array}{cc}
	X_{2}^{a}\frac{\partial f^{A}}{\partial u^{a}}= 0& \qquad A\neq6\\
	\\
	X_{2}^{a}\frac{\partial f^{6}}{\partial u^{a}}=1& 
\end{array}
\right..
\end{equation*}
We observe that the first equations in I) state that $(X_{1}^{a})_{a\in \{1,\ldots,5\}}$ is a non zero solution of a homogeneous system whose matrix has to be degenerate, i.e. $D^{5}=0$. Analogously from II) $D^{6}=0$ follows.
Vice versa, if we suppose that $D^{5}=0$, then I) admits non trivial solutions and  $\xi_{1}=X_{1}^{a}\frac{\partial}{\partial u^{a}}$. Analogously, supposing $D^{6}=0$, from II) we obtain $\xi_{2}=X_{2}^{a}\frac{\partial}{\partial u^{a}}$.
 \end{proof}

The following is an immediate consequence of Proposition \ref{det1} and Proposition \ref{det2}. 
\begin{proposition}\label{condlight}
Let $M$ be a hypersurface of $\bar M$. $M$ is a characteristic lightlike hypersurface if and only if 
$\;-(D^{1})^{2}+(D^{2})^{2}-(D^{3})^{2}+(D^{4})^{2}=0$ and $D^{5}=D^{6}=0$.
\end{proposition}
Now, we describe a lightlike hypersurface $M$ in $\bar M$ such that $\xi_{1}$, $\xi_{2}$ and $\varphi(E)$ belong to $S(TM)$.
We consider
\[
	f(u^{1},u^{2},u^{3},u^{4},u^{5})=(u^{1}+u^{5},u^{2},u^{3},u^{1}+u^{5},u^{4},u^{5}).
\]
It easy to check that this map describes a hypersurface in $\bar M$ and
 $D^{1}=D^{4}=1, \quad D^{2}=D^{3}=D^{5}=D^{6}=0$,
hence such a hypersurface verifies Proposition \ref{condlight}. Then $M$ is a lightlike hypersurface and $\xi_{1}$ and $\xi_{2}$ belong to $\Gamma(TM)$. The natural basis of $TM$ is given by 
\begin{align*}
U_{1}&=\frac{\partial}{\partial u^{1}}=\frac{\partial}{\partial x^{1}}+\frac{\partial}{\partial y^{2}},\quad 
U_{2}=\frac{\partial}{\partial u^{2}}=\frac{\partial}{\partial x^{2}},\quad U_{3}=\frac{\partial}{\partial u^{3}}=\frac{\partial}{\partial y^{1}},\\ 
U_{4}&=\frac{\partial}{\partial u^{4}}=\frac{\partial}{\partial z^{1}}=\xi_{1},\quad U_{5}=
\frac{\partial}{\partial u^{5}}=\frac{\partial}{\partial x^{1}}+\frac{\partial}{\partial y^{2}}+\xi_{2}.
\end{align*} 
We deduce that $\xi_{2}=U_{5}-U_{1}$. Considering the vector field of $\bar M$  
\[
E=-\frac{\partial}{\partial x^{1}}-\frac{\partial}{\partial y^{2}}+y^{1}\xi_{1}+y^{1}\xi_{2},
\]
it is easy to check that $E$ belongs to $Rad(TM)$. In fact we have
$g(E,U_{a})=0$ for all $a\in\{1,\ldots,5\}$ and $ \quad g(E,E)=0$. We construct $N\in\Gamma(T\bar M)$ such that
$g(N,N)=0$, and $ g(N,E)=1$.
To this aim, we consider $Z=\frac{\partial}{\partial x^{1}}-y^{1}\xi_{1}-y^{1}\xi_{2}$, we get $g(Z,E)=\frac{1}{2}\neq0$, and putting
\[
N=\frac{1}{g(Z,E)}\left(Z-\frac{g(Z,Z)}{2g(Z,E)}E\right),
\]
we have $N=\frac{\partial}{\partial x^{1}}-\frac{\partial}{\partial y^{2}}-y^{1}\xi_{1}-y^{1}\xi_{2}$.

Now, to give a basis of the screen distribution, we compute $\varphi(E)$ and $\varphi(N)$ and we have
\[
\varphi(E)=-\frac{\partial}{\partial x^{2}}+\frac{\partial}{\partial y^{1}}-y^{2}\xi_{1}-y^{2}\xi_{2},\quad
\varphi(N)=-\frac{\partial}{\partial x^{2}}-\frac{\partial}{\partial y^{1}}-y^{2}\xi_{1}-y^{2}\xi_{2}.
\]
Hence $\varphi(E)=-U_{2}+U_{3}-y^{2}U_{4}-y^{2}(U_{5}-U_{1})$ and $\varphi(N)=-U_{2}-U_{3}-y^{2}U_{4}-y^{2}(U_{5}-U_{1})$ are linear combinations of $U_{\alpha}$, so they are tangent to $M$.
Thus, we define the characteristic screen distribution as
$S(TM)=<\xi_{1},\xi_{2},\varphi(E),\varphi(N)>$,
and it is easy to check that the restriction of $\bar{g}$ (or $g$) to $S(TM)$ 
 has index one and therefore it is a metric tensor field on $S(TM)$. 
Moreover, we find
\[S(TM)^{\bot}=<\frac{\partial}{\partial x^{1}}-y^{1}\xi_{1}-y^{1}\xi_{2},\; \frac{\partial}{\partial y^{2}}>.\]
Finally, being $Z\in\Gamma(S(TM)^{\bot})$, we deduce that $N$ is the vector field provided by Theorem \ref{teorsezioneN}.

Now, for the described lightlike hypersurface of $\bar M$ we compute its geometrical object and write the Gauss and Weingarten equations for it and for the screen distribution.

Let $\bar{\nabla}$ be the Levi-Civita connection of $\bar M$ and $\{E,\xi_{1},\xi_{2},U=-\varphi E,V=-\varphi N\}$ a basis of $TM$. Then, reminding that $B(X,E)=0$ for any $X\in\Gamma(TM)$ and that $B$ is symmetric, we compute $B$ with respect to the aforesaid basis and we get, being $\varepsilon_1 = \varepsilon_2 =1$,   
\begin{align*}
	B(X,\xi_{\alpha})&=\bar{g}(\bar{\nabla}_{X}\xi_{\alpha},E)=-\bar{g}(\varphi X,E)=\left\{
	\begin{array}{ll}
	0 & X=\xi_{\beta}, \qquad \beta\in\{1,2\}\\
	\bar{g}(\varphi^{2}E,E)=0 &X=U\\
	\bar{g}(\varphi^{2}N,E)=-1 &X=V
	\end{array}
	\right..
\end{align*}
Hence, $M$ is neither totally geodesic nor totally umbilical since $\bar g (V, \bar\xi_{\alpha})= -\bar g (\bar\varphi N,\bar\xi_{\alpha})=0$.

We compute the Christoffel symbols of the Levi-Civita connection $\bar{\nabla}$ in order to find $B(U,U)$ and $B(V,V)$. It is easy to check 
\begin{align*}
	\Gamma_{11}^{3}&=4y^{1},\quad  \Gamma_{12}^{3}=\Gamma_{21}^{3}=-2y^{2}, \quad \Gamma_{12}^{4}=\Gamma_{21}^{4}=2y^{1}, \quad \Gamma_{22}^{4}=-4y^{2}, \\
	 \Gamma_{13}^{5}&=\Gamma_{13}^{6}=\Gamma_{31}^{5}=\Gamma_{31}^{6}=\frac{1}{2}+2(y^{1})^{2},\quad
	 \Gamma_{13}^{1}=\Gamma_{31}^{1}=\Gamma_{14}^{2}=\Gamma_{41}^{2}=-2y^{1},\\
	 \Gamma_{23}^{1}&=\Gamma_{32}^{1}=\Gamma_{24}^{2}=\Gamma_{42}^{2}=2y^{2}, \quad \Gamma_{24}^{5}=\Gamma_{42}^{5}=\Gamma_{24}^{6}=\Gamma_{42}^{6}=-\frac{1}{2}+2(y^{2})^{2},\\	\Gamma_{14}^{5}&=\Gamma_{14}^{6}=\Gamma_{41}^{5}=\Gamma_{41}^{6}=\Gamma_{23}^{5}=\Gamma_{23}^{6}=\Gamma_{32}^{5}=\Gamma_{32}^{6}=-2y^{1}y^{2},\\	 \Gamma_{15}^{3}&=\Gamma_{51}^{3}=\Gamma_{16}^{3}=\Gamma_{61}^{3}=\Gamma_{25}^{4}=\Gamma_{52}^{4}=\Gamma_{26}^{4}=\Gamma_{62}^{4}=1,\\
	\Gamma_{35}^{5}&=\Gamma_{53}^{5}=\Gamma_{36}^{5}=\Gamma_{63}^{5}=
	\Gamma_{35}^{6}=\Gamma_{53}^{6}=\Gamma_{36}^{6}=\Gamma_{63}^{6}=y^{1},\\	\Gamma_{45}^{2}&=\Gamma_{54}^{2}=\Gamma_{46}^{2}=\Gamma_{64}^{2}=\Gamma_{35}^{1}=\Gamma_{53}^{1}=\Gamma_{36}^{1}=\Gamma_{63}^{1}=-1,\\
	\Gamma_{45}^{5}&=\Gamma_{54}^{5}=\Gamma_{46}^{5}=\Gamma_{64}^{5}=
	\Gamma_{45}^{6}=\Gamma_{54}^{6}=\Gamma_{46}^{6}=\Gamma_{64}^{6}=-y^{2},
\end{align*}
while the other $\Gamma_{ij}^{h}$ vanish.
Being $U=\frac{\partial}{\partial x^{2}}-\frac{\partial}{\partial y^{1}}+y^{2}(\frac{\partial}{\partial z^{1}}+\frac{\partial}{\partial z^{2}})$, we have $\bar{\nabla}_{U}U=0$,
hence, we deduce
$B(U,U)=\bar{g}(\bar{\nabla}_{U}U,E)=0$.
 Finally, being $V=\frac{\partial}{\partial x^{2}}+\frac{\partial}{\partial y^{1}}+y^{2}(\frac{\partial}{\partial z^{1}}+\frac{\partial}{\partial z^{2}})$, we obtain $\bar{\nabla}_{V}V=0$ and $B(V,V)=\bar{g}(\bar{\nabla}_{V}V,E)=0$ and  $M$ is  minimal.
 
Using (\ref{equ12}), we compute $\tau$. To this aim, we find the value of $\bar{\nabla}_{X}N$ for any $X\in\Gamma(TM)$. For any $\alpha\in\{1,2\}$, it is easy to check $\bar{\nabla}_{E}N=0$, $\bar{\nabla}_{\xi_{\alpha}}N=\frac{\partial}{\partial y^{1}}+\frac{\partial}{\partial x^{2}}+y^{2}(\frac{\partial}{\partial z^{1}}+\frac{\partial}{\partial z^{2}})$, $\bar{\nabla}_{U}N=\frac{\partial}{\partial z^{1}}+\frac{\partial}{\partial z^{2}}$ and $\bar{\nabla}_{V}N=0$.
Immediately we find
\begin{align*}
	\tau(\xi_{\alpha})&=\bar{g}(\bar{\nabla}_{\xi_{\alpha}}N,E)=0, \text{ for any } \alpha\in\{1,2\}\\
	\tau(E)&=\bar{g}(\bar{\nabla}_{E}N,E)=0,\quad\tau(U)=\bar{g}(\bar{\nabla}_{U}N,E)=0,\quad \tau(V)=\bar{g}(\bar{\nabla}_{V}N,E)=0.
\end{align*}
So, from (\ref{equ12}), we have $A_{N}X=-\bar{\nabla}_{X}N$,
for any $X\in\Gamma(TM)$.

Again using (\ref{equ12}), we compute $\nabla E$ and we get $\nabla_{\xi_{\alpha}}E=-\frac{\partial}{\partial y^{1}}+\frac{\partial}{\partial x^{2}}+y^{2}(\frac{\partial}{\partial z^{1}}+\frac{\partial}{\partial z^{2}})$, $\nabla_{E}E=0$, $\nabla_{U}E=0$ and $\nabla_{V}E=\frac{\partial}{\partial z^{1}}+\frac{\partial}{\partial z^{2}}$.
Using (\ref{equ21}) and being $\tau(X)=0$ for any $X\in\Gamma(TM)$, we have
$\overset{*} {A}_{E}X=-\nabla_{X}E$, for any $X\in\Gamma(TM)$.
Finally, recalling that $S(TM)=<\xi_{1},\xi_{2},U,V>$ and $C(X,PY)=\bar{g}(\overset{*}{h}(X,PY),N)$, we value $C$ and we find
\[
C(X,PY)=\bar{g}(\overset{*}{h}(X,PY),N)=\bar{g}(\nabla_{X}PY,N)=\bar{g}(\bar{\nabla}_{X}PY,N).
\]
For any $\alpha\in\{1,2\}$, we have
\begin{align}\label{CperD0}
	C(X,\xi_{\alpha})&=\bar{g}(\bar{\nabla}_{X}\xi_{\alpha},N)=-\bar{g}(\varphi X,N)=\bar{g}(X,\varphi N)
	=\left\{
	\begin{array}{ccc}
	0&\quad \text{for }X=\xi_{\beta},&\,\beta\in\{1,2\}\\
	0&\quad \text{for }X=V\\
	-1&\quad \text{for }X=U\\
	0&\quad \text{for }X=E
	\end{array}
	\right..
\end{align}
Being $\bar{\nabla}_{\xi_{\alpha}}V=\frac{\partial}{\partial y^{2}}-\frac{\partial}{\partial x^{1}}+y^{1}(\frac{\partial}{\partial z^{1}}+\frac{\partial}{\partial z^{2}})$, we obtain
$C(\xi_{\alpha},V)=\bar{g}(\bar{\nabla}_{\xi_{\alpha}}V,N)=0$.\\
Analogously, $\bar{\nabla}_{U}V=0$ and $\bar{\nabla}_{V}V=0$ imply
$C(U,V)=\bar{g}(\bar{\nabla}_{U}V,N)=0$, $C(V,V)=\bar{g}(\bar{\nabla}_{V}V,N)=0$.
As $\bar{\nabla}_{E}V=-\frac{\partial}{\partial z^{1}}-\frac{\partial}{\partial z^{2}}$, we get
$C(E,V)=\bar{g}(\bar{\nabla}_{E}V,N)=0$.
Since $\bar{\nabla}_{V}U=0$ and $\bar{\nabla}_{U}U=0$, we have
$C(V,U)=\bar{g}(\bar{\nabla}_{V}U,N)=0, \quad C(U,U)=\bar{g}(\bar{\nabla}_{U}U,N)=0$.
Finally, $\bar{\nabla}_{E}U=0$ and $\bar{\nabla}_{\xi_{\alpha}}U=\frac{\partial}{\partial y^{2}}+\frac{\partial}{\partial x^{1}}-y^{1}(\frac{\partial}{\partial z^{1}}+\frac{\partial}{\partial z^{2}})$ for any $\alpha\in\{1,2\}$, we find
\[
C(E,U)=\bar{g}(\bar{\nabla}_{E}U,N)=0, \quad C(\xi_{\alpha},U)=\bar{g}(\bar{\nabla}_{\xi_{\alpha}}U,N)=-1.
\]

We check some properties of $D_{0}$. To this aim, we note that the distribution $D_{0}=<\xi_{1},\xi_{2}>$ coincides with $\ker(\bar \varphi)$ so it is integrable, totally geodesic and flat with respect to $\bar\nabla$ (\cite{LP}). It is also easy to check that $D_{0}$ verifies the conditions of Theorem \ref{integrD0}.
 Moreover, from (\ref{2°forma}), we obtain
\begin{align*}
\overset{\circ}{h}(\xi_{\alpha},\xi_{\beta}) =-C(\xi_{\alpha},\varphi\xi_{\beta})
\varphi E-B(\xi_{\alpha},\varphi\xi_{\beta})\varphi N +C(\xi_{\alpha},\xi_{\beta})E=0.
\end{align*}
and hence $\overset{\circ}{h}|_{D_{0}\times D_{0}}$ vanishes, while the second fundamental form $\overset{\circ}{h} :\Gamma(TM)\times\Gamma(D_{0})\rightarrow\Gamma(\mathcal{F})$ does not vanish, for example $\overset{\circ}{h}(E,\xi_{\alpha})=-\varphi E$, for any $\alpha\in\{1,2\}$. Furthermore, for any $\alpha\in\{1,2\}$, $\overset{\circ}{h}(\xi_{\alpha},\xi_{\alpha})=0$, then we have 
	$trace(\overset{\circ}{h})=0$
	i.e. $D_{0}$ is a minimal distribution in $M$, with respect to $\nabla$.

\medskip
Finally, we consider the indefinite $\mathcal{S}$-manifold $\bar M =({\mathbb R}^6_2,\varphi,\xi_{1},\xi_{2},\eta^{1},\eta^{2},g)$, described in Example 4.1 in \cite{LP} and given by
 $\xi_{\alpha}=\frac{\partial}{\partial z^{\alpha}}$, $\eta^{\alpha}=dz^{\alpha}-\sum_{i=1}^{2}y^{i}dx^{i}$, $\alpha\in\{1,2\}$, $ g=-\sum_{\alpha=1}^{2}\eta^{\alpha}\otimes\eta^{\alpha}+\frac{1}{2}\sum_{i=1}^{2}((dx^{i})^{2}+(dy^{i})^{2})$,
and $\varphi$ given by the matrix
\begin{equation*}
F=
\left(
\begin{array}{ccc}
0    & I_{2}     & 0 \\
-I_{2} & 0  & 0\\
0  & Y  & 0
\end{array}
\right),\quad {\rm where}\quad
Y=
\left(
\begin{array}{cc}
y^{1} & y^{2} \\
 y^{1} & y^{2}
\end{array}
\right).
\end{equation*}
In this case, the Killing characteristic vector fields are both timelike and we obtain
\begin{align*}
M_{11}&=M_{22}=M_{33}=M_{44}=\frac{1}{8},\quad M_{55}=M_{66}=-\frac{1}{16}+\frac{1}{8}(y^{1})^{2}+\frac{1}{8}(y^{1})^{2},
\\
M_{15}&=\frac{1}{8}y^{1},\quad M_{16}=-\frac{1}{8}y^{1},\quad M_{25}=-\frac{1}{8}y^{2},\quad M_{26}=\frac{1}{8}y^{2},
\quad 
M_{56}=-\frac{1}{8}(y^{1})^{2}-\frac{1}{8}(y^{2})^{2},
\end{align*}
while the other $M_{AB}$ vanish.
Using (\ref{prop1}) and Proposition \ref{prop2}, we can prove the following proposition omitting its proof because it is analogous to that of Proposition \ref{det1}.
\begin{proposition} 
Let $M$ be a hypersurface of $\bar M$. Then, $M$ is a lightlike hypersurface if and only if
\begin{align*}	\frac{1}{2}(&(D^{1})^{2}+(D^{2})^{2}+(D^{3})^{2}+(D^{4})^{2}+(-\frac{1}{2}+(y^{2})^{2}+(y^{1})^{2})((D^{5})^{2}+(D^{6})^{2}))\\
	& -y^{2}(D^{2}D^{5}-D^{2}D^{6})+((y^{1})^{2}+(y^{2})^{2})D^{5}D^{6}+y^{1}(D^{1}D^{5}-D^{1}D^{6})=0.
\end{align*}
\end{proposition}
Now, we note that the proof of Proposition \ref{det2} does not depend on the metric, therefore it also holds in this case. Thus, omitting the proof, we have
\begin{proposition}
Let $M$ be a hypersurface of $\bar M$. $M$ is a characteristic lightlike hypersurface if and only if 
$\;(D^{1})^{2}+(D^{2})^{2}+(D^{3})^{2}+(D^{4})^{2}=0$ and $D^{5}=D^{6}=0$. Hence there are not characteristic lightlike hypersurfaces in $\bar M$.
\end{proposition}
 
\medskip
As regards the Example 4.3 in \cite{LP} which describes a $4$-dimensional Lorentzian $\mathcal{S}$-manifold with two characteristic vector fields of different causal type
Remark \ref{dimensione} ensures that there are not characteristic lightlike hypersurfaces, since $n=1$.

Authors address:\\
Department of Mathematics,  University of Bari \\
  Via E. Orabona 4,\\
  I-70125 Bari (Italy) \\
  {\texttt{brunetti@dm.uniba.it};\;  \texttt{pastore@dm.uniba.it}}

\end{document}